\documentclass[12pt]{article}
\usepackage{latexsym}
\usepackage{epsfig}
\usepackage{amsfonts,amssymb,amsmath}

\title{A Ferguson - Klass - LePage series representation of\\
multistable multifractional motions and related processes}
\author{R. Le Gu\'evel\\
\small{{\em Universit\'{e} Nantes, Laboratoire de Math\'{e}matiques
Jean Leray UMR CNRS 6629}}\\
\small{{\em 2 Rue de la Houssini\`{e}re - BP 92208 - F-44322 Nantes Cedex 3, France}}\\
\small{{\em ronan.leguevel@univ-nantes.fr}}\\
\small{ and } \\
J. L\'{e}vy V\'{e}hel \\
\small{{\em Projet APIS, INRIA Saclay, Parc Orsay Universit\'e}}\\
\small{{\em 4 rue Jacques Monod - Bat P - 91893 Orsay Cedex, France}}\\
\small{{\em jacques.levy-vehel@inria.fr}}}

\date{}
\def\bbbr{{\bf R}} 

\newtheorem{theo}{Theorem}

\newtheorem{prop}[theo]{Proposition}

\newtheorem{lem}[theo]{Lemma}
\newtheorem{cor}[theo]{Corollary}

\newcommand\E{\mbox{\sf E}}

\newcommand\e{\eta}

\newcommand\efd{\stackrel{fdd}{=}}

\newcommand\ed{\stackrel{d}{=}}
\renewcommand\P{{\sf P}}

\newcommand{\one}{\ifmmode {\sf 1}\hspace{-.26em}{\sf
l}\hspace{-.35em}{\sf \_} \else ${\sf 1}\hspace{-.26em}{\sf
l}\hspace{-.35em}{\sf \_}$ \fi}

\newcommand\varep{\varepsilon}

\newcommand\ho{H\"{o}lder }

\renewcommand{\Box}{\mbox{\rule{1ex}{1ex}}}

\begin{document}

\maketitle

\begin{abstract}
\noindent The study of non-stationary processes whose local form has controlled
properties is a fruitful and important area of research, both in theory and applications.
In \cite{FLV}, a particular way of constructing such processes was investigated,
leading in particular to {\it multifractional multistable processes}, which were built
using sums over Poisson processes. We present here a different construction
of these processes, based on the Ferguson - Klass - LePage
series representation of stable processes. We consider various particular
cases of interest, including multistable L\'evy motion, multistable reverse
Ornstein-Uhlenbeck process, log-fractional multistable motion and linear multistable
multifractional motion. We also show that the processes defined here have the
same finite dimensional distributions as the corresponding processes constructed
in \cite{FLV}. Finally, we display numerical experiments showing graphs of synthesized
paths of such processes.

\end{abstract}

\vspace{1cm}


\vspace{1cm}

\section{Introduction}

This work deals with a general method for building
stochastic processes for which certain aspects of the local form are prescribed.
We will mainly be interested here in local \ho regularity and local intensity
of jumps, but our construction allows in principle to control
other properties that could be of interest. Our approach
is in the same spirit as the one proposed in \cite{FLV}, but it
uses different methods. In particular, in \cite{FLV},
multistable processes, that is localisable processes which
are locally $\alpha$-stable, but where the index of stability
$\alpha$ varies with time, were constructed using sums over
Poisson processes. We present here an alternative construction
of such processes, based on the Ferguson - Klass - LePage
series representation of stable stochastic processes \cite{FK,LP1,LP2}.
This representation is a powerful tool for the study of various
aspects of stable processes, see for instance \cite{BJP,JR}.
A comprehensive reference for the properties of this representation
that will be  needed here is \cite{ST}.

\medskip

Stochastic processes where the local H\"{o}lder regularity
varies with a parameter $t$ are interesting
both from a theoretical and practical point of view.
A well-known example is multifractional
Brownian motion (mBm), where the Hurst index $h$ of fractional
Brownian motion \cite{K40,MV} is replaced by a functional
parameter $h(t)$, permitting the H\"{o}lder exponent to
vary in a prescribed manner \cite{AL2,BJR,EH,PL}.
This allows in addition local regularity
and long range dependence to be decoupled
to give sample paths that are both highly
irregular and highly correlated, a useful feature
for instance in terrain or TCP traffic modeling.

However, local regularity, as measured by the \ho exponent, is not
the only local feature of a process that is useful
in theory and applications. Jump characteristics also
need to be accounted for, {\it e.g.} for studying processes
with paths in $D(\bbbr)$ (the space of c\`{a}dl\`{a}g functions,
{\it i.e.} functions which are continuous on the right and
have left limits at all $t \in T$). This has applications for
instance in the modeling of financial or medical data.
Stable non-Gaussian processes yield relevant models for
data containing discontinuities, with the stability
index $\alpha$ controlling the distribution of jumps.

Just for the same reason why it is interesting to consider stochastic processes whose
local \ho exponent changes in a controlled manner, tractable models where
the ``jump intensity'' $\alpha$ is allowed to vary in time
are needed, for instance to obtain a more accurate description of
some aspects of the local structure of functions in $D(\bbbr)$.

The approach described in this work allows in particular
to construct processes where both $h$ and $\alpha$
evolve in time in a prescribed way. Having two functional parameters allows
to finely tune the local properties of these processes. This may prove
useful to model two distinct aspects of financial risk,
to describe epileptic episodes in EEG
where for some periods there may be
only small jumps and at other instants very large ones,
or to study textured images where both H\"{o}lder regularity and
the distribution of discontinuities  vary.

\medskip

Let us now recall the definition of
a localisable process \cite{Fal5,Fal6}: $Y =\{Y(t): t \in \bbbr\}$
is said to be $h-$localisable at $u$ if there exists an $h \in \bbbr$
and a non-trivial limiting process $Y_{u}'$ such that
\begin{equation}
\lim_{r \to 0}\frac{Y(u+rt) -Y(u)}{r^{h}} = Y_{u}'(t).
\label{locform1}
\end{equation}
(Note $Y_{u}'$ may and in general will vary with $u$.) When the limit exits,
$Y_{u}'=\{Y_{u}'(t): t\in \bbbr\}$ is termed the {\it local form} or tangent
process of $Y$ at $u$ (see \cite{BJR,PL} for similar notions).
The limit (\ref{locform1}) may
be taken in several ways. In this work, we will only
deal with the case where convergence occurs
in finite dimensional distributions (equality in
finite dimensional distributions is denoted $\efd$). When convergence
takes place in distribution, the process is called
{\it strongly $h$-localisable} (equality in
distributions is denoted $\ed$).

As mentioned above, a now classical example is
multifractional Brownian motion $Y$
which ``looks like'' index-$h(u)$ fractional
Brownian motion close to time $u$ but where $h(u)$ varies, that is
\begin{equation}
\lim_{r \to 0}\frac{Y(u+rt) -Y(u)}{r^{h}} = B_{h(u)}(t) \label{exfbm}
\end{equation}
where $B_{h}$ is index-$h$ fractional Brownian motion.
A generalization of mBm, where the Gaussian
measure is replaced by an $\alpha$-stable one, leads to
multifractional stable processes, where the
local form is an $h(u)$-self-similar linear
$\alpha$-stable motion \cite{ST1,ST2}.

The $h$-local form $Y_{u}'$ at $u$, if it exists,
must be $h$-self-similar, that is $Y_{u}'(rt)\ed r^{h}Y_{u}'(t)$
for $r>0$. In addition, as shown in \cite{Fal5,Fal6}, under quite general
conditions, $Y_{u}'$ must also have stationary increments at almost all $u$ at
which there is strong localisability.  Thus, typical local forms
are self-similar with stationary increments (sssi), that is
$r^{-h}(Y_{u}'(u+rt) -Y_{u}'(u)) \ed Y_{u}'(t)$ for
all $u$ and $r>0$. Conversely, all sssi processes are localisable.
Classes of known sssi processes include fractional Brownian
motion, linear fractional stable motion and $\alpha$-stable L\'{e}vy
motion, see \cite{EM,ST}.

Similarly to \cite{FLV}, our method for constructing localisable processes
is to make use of stochastic fields
$\{X(t,v), (t,v)\in \bbbr^2\}$, where $t$ is time, and where the process
$ t \mapsto X(t,v)$ is localisable for all $v$. This field will allow to control
the local form of a `diagonal' process
$Y=\{X(t,t): t\in \bbbr\}$. For instance, in the case
of mBm, $X$ will be a field of fractional Brownian motions, {\it i.e.}
$X(t,v) = B_{h(v)}(t)$, where $h$ is a smooth function of $v$
ranging in $[a,b] \subset (0,1)$. This is the approach that was used originally
in \cite{AL2} for studying mBm. From a heuristic point of view, taking
the diagonal of such a stochastic field constructs a new process  with
local form depending on $t$ by piecing together
known localisable processes. In other words, we shall use random fields
$\{X(t,v):(t,v) \in \bbbr^2\}$ such that
for each $v$ the local form $X_{v}'(\cdot,v)$ of $X(\cdot,v)$ at $v$
is the desired local
form $Y_{v}'$ of $Y$ at $v$.  An easy situation is when, for each $v$, the
process $\{X(t,v): t\in \bbbr\}$ is sssi, since this automatically entails
localisability.

It is clear that, in this approach, the structure of
$X(\cdot,v)$ for $v$ in a neighbourhood of $u$ will be crucial to
determine the local
behaviour of $Y$ near $u$. A simple way to control this structure
is to define the random field as
an integral or sum of functions that depend on $t$ and $v$ with
respect to a single underlying random measure so as to
provide the necessary correlations.

General criteria that guarantee the transference of
localisability from the $X(\cdot,v)$ to $Y=\{X(t,t): t \in \bbbr\}$
were obtained in \cite{FLV}. We will make use of the
following one:

\begin{theo}\label{theoflv}
Let $U$ be an interval with $u$ an interior point.
Suppose that for some $0<h<\eta$ the process  $\{X(t,u), t \in
U\}$
is $h$-localisable at $u\in U$
with local form $X_{u}'(\cdot,u)$ and
\begin{equation}
\P(|X(v,v) - X(v,u)|\geq |v-u|^{\e}) \to 0 \label{cond}
\end{equation}
as $v\to u$.  Then $Y=\{X(t,t) : t\in U\}$ is $h$-localisable at $u$
with $Y_{u}'(\cdot) = X_{u}'(\cdot,u)$.
\end{theo}

\medskip

In the sequel, we shall consider specific classes of random fields
and use Theorem \ref{theoflv} to build localisable processes with
interesting local properties. As a particular case, we will study
multifractional multistable processes, where both the local \ho regularity
and intensity of jumps will evolve in a controlled manner.

The remaining of this article is organized as follows: we first
collect some notations in section \ref{notations}.
We then build localisable processes using a series representation that
yields the necessary flexibility required for our purpose.
We need to distinguish between the
situations where the underlying space is finite (section \ref{finite}),
or merely $\sigma-$ finite (section \ref{sigfinite}). In each case,
we define a random field depending on a ``kernel'' $f$, and give
conditions on $f$ ensuring localisability of the diagonal process.
We then consider in section \ref{Examples} some examples:
multistable L\'evy motion, multistable reverse Ornstein-Uhlenbeck process,
log-fractional multistable motion and linear multistable multifractional motion.
Section \ref{fdd} is devoted to computing the finite dimensional
distributions of our processes, and proving that they are the
same as the ones of the corresponding processes constructed in \cite{FLV}.
Finally, section \ref{numexp} displays graphs of certain localisable
processes of interest, in particular multifractional multistable ones.

\medskip
Before we proceed, we note that constructing localisable processes
using a stochastic field composed of sssi processes is obviously
not the only approach that one can think of.
It is for instance possible to follow a rather different path and construct
localisable processes from moving average ones by imposing
conditions on the kernel defining the moving average. See \cite{FLGLV} for details.

\section{Notations}\label{notations}

We refer the reader to the first chapters of \cite{ST} for basic notions on stable random variables and stable processes. In particular, recall that a process $\{X(t):t\in T\}$, where $T$ is
generally a subinterval of $\bbbr$,
is called $\alpha$-{\it stable} $(0<\alpha \leq 2)$
if all its finite-dimensional distributions are
$\alpha$-stable. Many stable processes admit
a stochastic integral representation as follows.
Write $S_{\alpha}(\sigma,\beta,\mu)$ for the $\alpha$-stable distribution
with scale parameter $\sigma$, skewness $\beta$ and shift-parameter
$\mu$; we will assume throughout that $\mu=0$.  Let $(E,{\cal E},m)$ be a
sigma-finite measure space.
Taking $m$  as the control measure and $\beta: E \to [-1,1]$ a
measurable function, this defines an
$\alpha$-stable random measure $M$ on $E$ such that
for $A\in {\cal E}$ we have that
$M(A) \sim S_{\alpha}\left(m(A)^{1/\alpha},
\int_{A}\beta(x)m(dx)/m(A),0\right)$.  If $\beta =0$, the process is
termed {\it symmetric} $\alpha$-stable, or $S\alpha S$.

Let
$$    {\cal F}_{\alpha}\equiv {\cal F}_{\alpha}(E,{\cal E}, m)
= \{ f: f \mbox{ is measurable and } \|f\|_{\alpha} < \infty\},$$
where $\|\,\|_{\alpha}$ is the quasinorm (or norm if $1<\alpha \leq
2$) given by
\begin{equation}
\|f\|_{\alpha}
=\left\{
\begin{array}{cc}
    \left( \int_E |f(x)|^{\alpha}m(dx)\right)^{1/\alpha} & (\alpha\neq 1)  \\
     \int_E |f(x)|m(dx)
     + \int_E |f(x)\beta(x) \ln |f(x)| | m(dx)  & (\alpha = 1)
\end{array}
 \right.
    \label{normdef}
\end{equation}
The stochastic
integral of $f\in  {\cal F}_{\alpha}(E,{\cal E}, m)$ with respect to
$M$ then exists  \cite[Chapter 3]{ST} with
\begin{equation}
I(f) =\int_E f(x)M(dx)\sim
S_{\alpha}(\sigma_{f},\beta_{f},0),\label{alint}
\end{equation}
where
$$\sigma_{f}=\|f\|_{\alpha}, \quad
\beta_{f} = \frac{\int f(x)^{<\alpha>}\beta(x) m(dx)}
{\|f\|_{\alpha}^{\alpha}},$$
and $a^{<b>} \equiv \mbox{sign}({a}) |a|^{b}$, see \cite[Section
3.4]{ST}.
In particular,
\begin{equation}
\E|I(f)|^{p} = \left\{
\begin{array}{cc}
   c(\alpha,\beta,p)\|f\|_{\alpha}^{p}  &  (0<p<\alpha)  \\
   \infty  & (p \geq \alpha)
\end{array}\right. \label{alintexp}
\end{equation}
where $c(\alpha,\beta,p)<\infty$, see \cite[Property 1.2.17]{ST}.

In this work, we will only consider the symmetric case, {\it i.e.}
we take $\beta\equiv0$ for the remaining of the article. We believe
most results should have a counterpart in the
non-symmetric case, although the proofs would
probably have to be much more involved.

\section{A Ferguson - Klass - LePage series representation of localisable processes in the finite measure space case}\label{finite}

A well-known representation of stable random variables is the
Ferguson - Klass - LePage series one \cite{BJP,FK,LP1,LP2,JR}.
This representation is particularly adapted for our purpose since,
as we shall see, it allows for
easy generalization to the case of varying $\alpha$.

In this work, we will use the following version:

\begin{theo} (\cite[Theorem 3.10.1]{ST}) \label{tst}

Let $(E,{\cal E},m)$ be a finite measure space, and $M$ be a symmetric $\alpha$-stable random measure with $\alpha \in (0,2)$ and finite control measure $m$. Let $(\Gamma_i)_{i \geq 1}$ be a sequence of arrival times of a Poisson process with unit arrival time, $(V_i)_{i \geq 1}$ be a sequence of i.i.d. random variables with distribution $\hat m = m/m(E)$ on $E$, and $(\gamma_i)_{i \geq 1}$ be a sequence of i.i.d. random variables with distribution $P(\gamma_i=1)=P(\gamma_i=-1)=1/2$. Assume finally that the three sequences $(\Gamma_i)_{i \geq 1}$, $(V_i)_{i \geq 1}$, and $(\gamma_i)_{i \geq 1}$ are independent. Then, for any $f \in {\cal F}_{\alpha}(E,{\cal E}, m)$,
\begin{equation}\label{serep}
\int_E f(x)M(dx) \ed \left( C_\alpha m(E) \right)^{1/\alpha} \sum_{i=1}^{\infty} \gamma_i \Gamma_i^{-1/\alpha}f(V_i),
\end{equation}
\end{theo}

where $C_\alpha := \frac{1-\alpha}{\Gamma(2-\alpha)\cos(\pi \alpha/2)}$ for $\alpha \neq 1$, $C_1 = 2/\pi$ (Theorem 3.10.1 in \cite{ST} is more general, as it extends to non-symmetric stable processes, that are not considered here.) As mentioned above, a relevant feature of this representation for us is that the distributions of all random variables appearing in the sum are independent of $\alpha$. We will use (\ref{serep}) to contruct processes with varying $\alpha$ as described in the following Theorem.
\begin{theo}\label{msspfm}
Let $(E,{\cal E},m)$ be a finite measure space. Let $\alpha$ be a $C^1$ function defined on $\bbbr$ and ranging in $[c,d] \subset (0,2)$. Let $b$ be a $C^1$ function defined on $\bbbr$. Let $f(t,u,.)$ be a family of functions such that, for all $(t,u) \in \bbbr^2$, $f(t,u,.) \in {\cal F}_{\alpha(u)}(E,{\cal E}, m)$. Let $(\Gamma_i)_{i \geq 1}$ be a sequence of arrival times of a Poisson process with unit arrival time, $(V_i)_{i \geq 1}$ be a sequence of i.i.d. random variables with distribution $\hat m = m/m(E)$ on $E$, and $(\gamma_i)_{i \geq 1}$ be a sequence of i.i.d. random variables with distribution $P(\gamma_i=1)=P(\gamma_i=-1)=1/2$. Assume finally that the three sequences $(\Gamma_i)_{i \geq 1}$, $(V_i)_{i \geq 1}$, and $(\gamma_i)_{i \geq 1}$ are independent.
Consider the following random field:

\begin{equation}\label{msfm}
X(t,u)= b(u)(m(E))^{1/\alpha(u)}C^{1/\alpha(u)}_{\alpha(u)} \sum_{i=1}^{\infty} \gamma_i \Gamma_i^{-1/\alpha(u)} f(t,u,V_i),
\end{equation}
where $C_{\alpha} = \left( \int_{0}^{\infty} x^{-\alpha} \sin (x)dx \right)^{-1}$.
Assume that $X(t,u)$ (as a process in $t$) is localisable at $u$ with exponent $h \in (0,1)$ and local form $X'_u(t,u)$. Assume in addition that:

\begin{itemize}
    \item (C1) The family of functions $v \to f(t,v,x)$ is differentiable for all $(v,t)$ in a neighbourhood of $u$ and almost all $x$ in $E$. The derivatives of $f$ with respect to $v$ are denoted by $f'_v$.
    \item (C2) There exists $\varep >0$ such that:
\begin{equation}\label{kercond1f}
\sup_{t \in B(u,\varep)}  \int_E \sup_{w \in B(u,\varep)} (|f(t,w,x)|^{\alpha(w)}) \hspace{0.1cm} \hat m(dx) < \infty.
\end{equation}
	\item (C3) There exists $\varep >0$ such that:
\begin{equation}
\sup_{t \in B(u,\varep)}  \int_E \sup_{w \in B(u,\varep)} (|f'_v(t,w,x)|^{\alpha(w)}) \hspace{0.1cm} \hat m(dx) < \infty.
\end{equation}
	\item (C4) There exists $\varep >0$ such that:
\begin{equation}\label{kercond2f}
\sup_{t \in B(u,\varep)}  \int_E \sup_{w \in B(u,\varep)} \left[\left|f(t,w,x)\log|f(t,w,x)|\right|^{\alpha(w)}\right] \hspace{0.1cm} \hat m(dx) < \infty.
\end{equation}

\end{itemize}

Then $Y(t) \equiv X(t,t)$ is localisable at $u$ with exponent $h$ and local form $Y'_u(t)=X'_u(t,u)$.

\end{theo}

\medskip

\noindent{\it Proof}

The function $u \mapsto C^{1/\alpha(u)}_{\alpha(u)}$ is $C^1$ since $\alpha(u)$ ranges in $[c,d] \subset (0,2)$. We shall denote $a(u) = b(u)(m(E))^{1/\alpha(u)} C^{1/\alpha(u)}_{\alpha(u)}$. The function $a$ is thus also $C^1$. We want to apply Theorem 1.1. With that in view, we estimate, for
$v \in B(u,\varepsilon)$ (the ball centered at $u$ with radius $\varepsilon$),

$$X(v,v)-X(v,u) =: \sum_{i=1}^{\infty} \gamma_i (\Phi_i(v)-\Phi_i(u)) + \sum_{i=1}^{\infty}\gamma_i (\Psi_i(v)-\Psi_i(u)),$$

where
 \begin{displaymath}
\Phi_i(u) = a(u) i^{-1/\alpha(u)} f(v,u,V_i)
\end{displaymath}
and
\begin{displaymath}
\Psi_i(u)=a(u) \left( \Gamma_i^{-1/\alpha(u)} - i^{-1/\alpha(u)} \right) f(v,u,V_i).
\end{displaymath}
The reason for introducing the $\Phi_i$ and the $\Psi_i$ is that the
random variables $\Gamma_i$ are not independent, which complicates their study.
We shall decompose the sum involving the $\Phi_i$ into series of independent
random variables which will be dealt with using the three series theorem.
The sum involving the $\Psi_i$ will be studied by taking advantage of the
fact that, for large enough $i$, each $\Gamma_i$ is ``close'' to $i$ in some
sense.

In the sequel, since only the values of $\alpha$ inside $B(u,\varepsilon)$
matter, we shall agree by convention that $c$ denotes in fact
$\inf_{v \in B(u,\varepsilon)} \alpha(v)$ and likewise
$d = \sup_{v \in B(u,\varepsilon)} \alpha(v)$. Note that, by decreasing
$\varepsilon$, $d-c$ may be made arbitrarily small.

Thanks to the assumptions on $a $ and $f$, $\Phi_i$ and $\Psi_i$ are differentiable and one computes:
$$\Phi_i'(u) = a'(u) i^{-1/\alpha(u)} f(v,u,V_i) + a(u)i^{-1/\alpha(u)} f'_u(v,u,V_i) + a(u)\frac{\alpha'(u)}{\alpha(u)^2} \log(i) i^{-1/\alpha(u)} f(v,u,V_i),$$
and
$$
\Psi_i'(u) = a'(u) \left( \Gamma_i^{-1/\alpha(u)} - i^{-1/\alpha(u)} \right) f(v,u,V_i) + a(u)\left( \Gamma_i^{-1/\alpha(u)} - i^{-1/\alpha(u)} \right) f'_u(v,u,V_i)
$$
$$+ a(u)\frac{\alpha'(u)}{\alpha(u)^2} \left( \log(\Gamma_i) \Gamma_i^{-1/\alpha(u)} - \log(i) i^{-1/\alpha(u)} \right) f(v,u,V_i).
$$

The mean value theorem yields that there exists a sequence of independent random numbers $w_i \in [u,v]$ (or $[v,u]$) and a sequence of random numbers $x_i \in [u,v]$ (or $[v,u]$) such that:
$$X(v,u)-X(v,v) = (u-v)\sum_{i=1}^{\infty} (Z_i^1+Z_i^2+Z_i^3) + (u-v)\sum_{i=1}^{\infty} (Y_i^1+Y_i^2+Y_i^3),$$
where
$$Z_i^1 = \gamma_i a'(w_i)i^{-1/\alpha(w_i)} f(v,w_i,V_i),$$
$$Z_i^2 = \gamma_i a(w_i)i^{-1/\alpha(w_i)} f'_u(v,w_i,V_i),$$
$$Z_i^3 = \gamma_i a(w_i)\frac{\alpha'(w_i)}{\alpha(w_i)^2} \log(i) i^{-1/\alpha(w_i)} f(v,w_i,V_i),$$
$$Y_i^1 = \gamma_i a'(x_i)\left( \Gamma_i^{-1/\alpha(x_i)} - i^{-1/\alpha(x_i)} \right) f(v,x_i,V_i),$$
$$Y_i^2 = \gamma_i a(x_i)\left( \Gamma_i^{-1/\alpha(x_i)} - i^{-1/\alpha(x_i)} \right) f'_u(v,x_i,V_i),$$
$$Y_i^3 = \gamma_i a(x_i)\frac{\alpha'(x_i)}{\alpha(x_i)^2} \left( \log(\Gamma_i) \Gamma_i^{-1/\alpha(x_i)} - \log(i) i^{-1/\alpha(x_i)} \right) f(v,x_i,V_i).$$
Note that each $w_i$ depends on $a,f,\alpha,u,v,V_i$, but not on $\gamma_i$. This remark will be useful in the sequel.

The remainder of the proof is divided into four steps. The first step will apply the three-series theorem to show that each series $\sum_{i=1}^{\infty}\limits Z_i^j, j=1,2,3,$ converges almost surely. In the
second step, we will prove that $\sum_{i=1}^{\infty}\limits Y_i^j$ also
converges almost surely for $j=1,2,3$. In the third step
we will prove that condition (\ref{cond}) is verified by  $\sum_{i=1}^{\infty}\limits Z_i^j, j=1,2,3$. Finally, step four will prove
the same thing for $\sum_{i=1}^{\infty}\limits Y_i^j, j=1,2,3$.

\medskip

{\bf First step: almost sure convergence of $\sum_{i=1}^{\infty}\limits Z_i^j, j=1,2,3$.}

\smallskip

Consider $Z^1 = \sum_{i=1}^{\infty}\limits Z_i^1$. Fix $\lambda>0$. We shall deal successively
with the three series involved the three-series theorem.

\smallskip

\underline {First series}: $S_1=\sum_{i=1}^{\infty}\limits \P(|Z_i^1|>\lambda)$.

\begin{eqnarray*}
\P(|Z_i^1|>\lambda) & = & \P\left(|f(v,w_i,V_i)| > \frac{\lambda i^{1/\alpha(w_i)}}{|a'(w_i)|}\right)\\
& \leq &  \P\left(|f(v,w_i,V_i)|^{\alpha(w_i)} >  i \inf_{w \in
B(u,\varep)}
\left[\left(\frac{\lambda}{|a'(w)|}\right)^{\alpha(w)}\right]
\right).
\end{eqnarray*}
Note that, since $a'$ is bounded on the compact interval $[u,v]$,
$K :=\inf_{w \in B(u,\varep)}
\left[\left(\frac{\lambda}{|a'(w)|}\right)^{\alpha(w)}\right]$ is
strictly positive.

\begin{eqnarray*}
\P(|Z_i^1|>\lambda) & \leq & \P\left( \sup_{w \in B(u,\varep)}|f(v,w,V_i)|^{\alpha(w)} > Ki \right)\\
& = &  \P\left( \sup_{w \in B(u,\varep)}|f(v,w,V_1)|^{\alpha(w)} > Ki \right).\\
\end{eqnarray*}

Thus
\begin{eqnarray*}
\sum_{i=1}^{+\infty} \P(|Z_i^1|>\lambda) & \leq & \sum_{i=1}^{+\infty} \P\left( \sup_{w \in B(u,\varep)}|f(v,w,V_1)|^{\alpha(w)} > Ki \right)\\
& \leq & \frac{1}{K} \E\left[ \sup_{w \in B(u,\varep)}|f(v,w,V_1)|^{\alpha(w)} \right]\\
& \leq & \frac{1}{K} \sup_{t \in B(u,\varep)} \int_E \sup_{w \in B(u,\varep)} (|f(t,w,x)|^{\alpha(w)}) \hspace{0.1cm} \hat m(dx)\\
& < & +\infty\\
\end{eqnarray*}

by (C2).

\smallskip

\underline {Second series}: $S_2^n=\sum_{i=1}^{n}\limits \E(Z_i^1\one\{|Z_i^1|\leq\lambda\})$.

\begin{eqnarray*}
\E(Z_i^1\one\{|Z_i^1|\leq\lambda\})&  = & \E(\gamma_i a'(w_i)i^{-1/\alpha(w_i)} f(v,w_i,V_i)\one\{|a'(w_i)i^{-1/\alpha(w_i)} f(v,w_i,V_i)|\leq\lambda\})\\
& = & \E(\gamma_i) \E(a'(w_i)i^{-1/\alpha(w_i)} f(v,w_i,V_i)\one\{|a'(w_i)i^{-1/\alpha(w_i)} f(v,w_i,V_i)|\leq\lambda\})\\
& = & 0,
\end{eqnarray*}

where we have used the facts that $\gamma_i$ is independent of $w_i, V_i$ and $\E(\gamma_i) = 0$. As a consequence,
$ \lim_{n \rightarrow +\infty} S_2^n=0$.

\smallskip

\underline {Third series}: The final series we need to consider is $S_3=\sum_{i=1}^{\infty}\limits \E\left[(Z_i^1\one\{|Z_i^1|\leq 1\})^2\right]$. Take $\lambda = 1$\footnote{Recall that, in the three
series theorem, for the series $\sum_{i=1}^\infty X_i$ to converge
almost surely, it is necessary that, for all $\lambda >0$, the
three series $\sum_{i=1}^\infty \P(|X_i|> \lambda),
\sum_{i=1}^\infty \E(X_i \one(|X_i| \leq \lambda))$, and
$\sum_{i=1}^\infty \mbox{Var}(X_i\one(|X_i| \leq \lambda))$
converge, and it is sufficient that they converge for {\it one}
$\lambda >0$, see, e.g. \cite{Pet}, Theorem 6.1.}.

\smallskip
\smallskip

Let $\eta$ be such that $d < \eta < 2$.

\begin{eqnarray*}
(Z_i^1\one\{|Z_i^1|\leq 1\})^2 & \leq & |Z_i^1|^{\eta} \one\{|Z_i^1|\leq 1\}\\
\E\left[(Z_i^1\one\{|Z_i^1|\leq 1\})^2\right] & \leq & \E\left[|Z_i^1|^{\eta} \one\{|Z_i^1|\leq 1\}\right]\\
& = & \int_{0}^{+\infty} \P(|Z_i^1|^{\eta} \one\{|Z_i^1|\leq 1\} > x) dx\\
& = & \int_{0}^{1} \P(|Z_i^1|^{\eta} \one\{|Z_i^1|\leq 1\} > x) dx\\
& \leq & \int_{0}^{1} \P(|Z_i^1|^{\eta} > x) dx.\\
\end{eqnarray*}
Now, for all $x$ in $(0,1)$,
\begin{eqnarray*}
\P(|Z_i^1|^{\eta} > x) & = & \P(|Z_i^1| > x^{1/ \eta})\\
 & \leq & \P \left( |f(v,w_i,V_i)|^{\alpha(w_i)} > i \frac{x^{\frac{\alpha(w_i)}{\eta}}}{|a'(w_i)|^{\alpha(w_i)}} \right)\\
 & \leq & \P\left( \sup_{w \in B(u,\varep)}|f(v,w,V_i)|^{\alpha(w)} > K'ix^{\frac{d}{\eta}} \right)\\
 & = & \P\left( \sup_{w \in B(u,\varep)}|f(v,w,V_1)|^{\alpha(w)} > K'ix^{\frac{d}{\eta}} \right),\\
\end{eqnarray*}
where $K' :=\inf_{w \in B(u,\varep)}
\left[\left(\frac{1}{|a'(w)|}\right)^{\alpha(w)}\right]$ is
strictly positive. Thus,

\begin{eqnarray*}
S_3 & \leq & \int_{0}^{1} \sum_{i=1}^{\infty} \P\left( \sup_{w \in B(u,\varep)}|f(v,w,V_1)|^{\alpha(w)} > K'ix^{\frac{d}{\eta}} \right) dx\\
& \leq & \int_{0}^{1} \frac{1}{K'x^{\frac{d}{\eta}}} \E (\sup_{w \in B(u,\varep)}|f(v,w,V_1)|^{\alpha(w)}) dx\\
& \leq & \frac{1}{K'} \left( \sup_{t \in B(u,\varep)} \int_E \sup_{w \in B(u,\varep)} (|f(t,w,x)|^{\alpha(w)}) \hspace{0.1cm} \hat m(dx) \right) \left(\int_{0}^{1} \frac{dx}{x^{\frac{d}{\eta}}} \right)\\
& < & +\infty.\\
\end{eqnarray*}

\bigskip

The case of the $Z^2= \sum_{i=1}^{\infty}\limits Z_i^2$ is treated similarly, since the conditions required on $(a',f)$ in the proof above are also satisfied by $(a,f'_u)$.

\bigskip

Consider finally $Z^3 = \sum_{i=1}^{\infty}\limits Z_i^3$. Fix $\lambda>0$.

\smallskip

\underline {First series}: $S_1=\sum_{i=1}^{\infty}\limits \P(|Z_i^3|>\lambda)$.

\begin{eqnarray*}
\P(|Z_i^3|>\lambda) & =  &   \P\left(|f(v,w_i,V_i)| > \frac{\lambda \alpha(w_i)^2 i^{1/\alpha(w_i)}}{|a(w_i) \alpha'(w_i)|\log i}\right) \\
& \leq &  \P\left(|f(v,w_i,V_i)|^{\alpha(w_i)} > K'' \frac{i}{(\log i )^{\alpha(w_i)}}\right),
\end{eqnarray*}
 where
$K'' := \inf_{w \in B(u,\varep)}
\left[\left(\frac{\lambda \alpha(w)^2}{|a(w) \alpha'(w)|}\right)^{\alpha(w)}\right]$ is strictly positive by the assumptions on $a,\alpha$ and $\alpha'$. In the sequel, $K$ will always denote a finite positive constant, that may however change from line to line.

Let $g_i(x) = \frac{Kx}{( \log x)^{\alpha(w_i)}}$ for $x \geq 1$ and $i \in \mathbb{N}^*$. For $x$ large enough and for all $i$, $g_i$ is strictly increasing and $ \lim_{x \rightarrow +\infty} g_i(x)= + \infty$. For $z$ large enough (independently of $i$),
\begin{eqnarray*}
g_i(z( \log z )^{\alpha(w_i)}) & = & \frac{Kz(\log z)^{\alpha(w_i)}}{(\log z + \alpha(w_i)\log \log z )^{\alpha(w_i)}}\\
 & \geq & \frac{Kz}{2}.\\
\end{eqnarray*}

Let $A > e$ be such that: $\forall z \geq A, \forall i \in \mathbb{N}^* , g_i^{-1}(z) \leq \frac{z}{K}(\log z)^{\alpha(w_i)}$.

Let $U_i =|f(v,w_i,V_i)|^{\alpha(w_i)}$.

\begin{eqnarray*}
\P(|Z_i^3|>\lambda) & \leq & \P( U_i > K\frac{i}{(\log i)^{\alpha(w_i)}})\\
& \leq & \P(\left( \{ U_i < A\} \cap \{U_i > K\frac{i}{(\log i)^{\alpha(w_i)}} \}\right) \cup \left( \{ U_i \geq A\} \cap \{ g_i^{-1}(U_i) > i \} \right))\\
& \leq & \P(A > U_i > K\frac{i}{(\log i)^d}) + \P(\{ U_i \geq A\} \cap \{ U_i |\log U_i|^{\alpha(w_i)} > Ki \})\\
& \leq & \P(A > U_i > K\frac{i}{(\log i)^d}) + \P(U_i |\log U_i|^{\alpha(w_i)} > Ki )\\
& \leq & \P(A > U_i > K\frac{i}{(\log i)^d}) + \P( \sup_{w \in B(u,\varep)} [|f(v,w,V_1)\log|f(v,w,V_1)||^{\alpha(w)}]  > \frac{K}{d^d}i).\\
\end{eqnarray*}

On the one hand, $\P(A > U_i > K\frac{i}{(\log i)^d}) $ vanishes for i large. On the other hand, by (C3),

\begin{eqnarray*}
 \sum_i \P\left( \sup_{w \in B(u,\varep)} [|f(v,w,V_1)\log|f(v,w,V_1)||^{\alpha(w)}] > \frac{K}{d^d}i \right) & \leq &\\
\frac{d^d}{K} \E \left( \sup_{w \in B(u,\varep)} [|f(v,w,V_1)\log|f(v,w,V_1)||^{\alpha(w)}] \right) < +\infty & & \\
\end{eqnarray*}
and thus $S_1 < + \infty$.

\smallskip

\underline {Second series}: $S_2^n=\sum_{i=1}^{n}\limits \E(Z_i^3\one\{|Z_i^3|\leq\lambda\})$.

For the same reason as in the case of $Z_i^1$ ({\it i.e.} $\gamma_i$ is independent of $w_i, V_i$), $S_2^n = 0$, $\forall n$.
\smallskip

\underline {Third series}: $S_3=\sum_{i=1}^{\infty}\limits \E((Z_i^3\one\{|Z_i^3|\leq 1\})^2)$. Take again $\lambda = 1$.

\smallskip
\smallskip

Let $\eta > 0$ be such that $d < \eta < 2$, and $\mu > 0$ be such that $0 < \mu <1 - \frac{d}{\eta}$. Using the same line of reasoning as in the case of $Z_i^1$\label{page10},

\begin{eqnarray*}
 \E(Z_i^3\one\{|Z_i^3|\leq 1\})^2 & \leq & \int_{0}^{1} \P(|Z_i^3|^{\eta} > x) dx.\\
\end{eqnarray*}
We have, for all $x$ in $(0,1)$,
\begin{eqnarray*}
\P(|Z_i^3|^{\eta} > x) & = & \P(|Z_i^3| > x^{1/ \eta})\\
& \leq &  \P\left(|f(v,w_i,V_i)|^{\alpha(w_i)} > K x^{\frac{d}{\eta}}\frac{i}{(\log i )^{\alpha(w_i)}}\right),
\end{eqnarray*}
 where
$K := \inf_{w \in B(u,\varep)}
\left[\left(\frac{ \alpha(w)^2}{|a(w) \alpha'(w)|}\right)^{\alpha(w)}\right]$ is strictly positive by the assumptions on $a,\alpha$ and $\alpha'$.

\begin{eqnarray*}
\P(|Z_i^3|^{\eta} > x) & \leq & \P(A > U_i > K x^{\frac{d}{\eta}} \frac{i}{(\log i)^d}) + \P( \sup_{w \in B(u,\varep)} [|f(v,w,V_1)\log|f(v,w,V_1)||^{\alpha(w)}]  > \frac{Kx^{\frac{d}{\eta}} }{d^d}i).\\
\end{eqnarray*}
Let us deal with the first term of the sum in the right hand side of the above inequality.
\begin{eqnarray*}
\P(A > U_i > K x^{\frac{d}{\eta}} \frac{i}{(\log i)^d}) & \leq & \P(A > U_i > K x^{\frac{d}{\eta}} i^{1-\mu})\\
\end{eqnarray*}
for $i$ large enough, {\it i.e.} $i > i^*$, where $i^*$ depends only on $d$ and $\mu$. As a consequence
\label{page11}:

\begin{eqnarray*}
 \sum_{i=i*}^{+\infty} \P(A > U_i > K x^{\frac{d}{\eta}} \frac{i}{(\log i)^d}) & \leq & \sum_{i=i*}^{+\infty} \P(A > U_i > K x^{\frac{d}{\eta}} i^{1-\mu})\\
 & \leq & \left( \frac{A}{K x^{\frac{d}{\eta}}} \right)^{\frac{1}{1-\mu}} \\
 & \leq & \frac{K'}{x^{\frac{d}{\eta(1-\mu)}}}.
\end{eqnarray*}
But $\frac{d}{\eta(1-\mu)} < 1$  and thus

\begin{eqnarray*}
\int_{0}^{1} \sum_{i=1}^{+\infty} \P(A > U_i > K x^{\frac{d}{\eta}} \frac{i}{(\log i)^d}) dx & < & +\infty.\\
\end{eqnarray*}

Now for the second term in the sum:
\begin{eqnarray*}
\sum_{i=1}^{+\infty} \P( \sup_{w \in B(u,\varep)} [|f(v,w,V_1)\log|f(v,w,V_1)||^{\alpha(w)}]  > \frac{Kx^{\frac{d}{\eta}} }{d^d}i) &  &\\
\leq  \frac{d^d}{Kx^{\frac{d}{\eta}}} \E \left( \sup_{w \in B(u,\varep)} [|f(v,w,V_1)\log|f(v,w,V_1)||^{\alpha(w)}]\right) & &\\
  \leq  \frac{K'}{x^{\frac{d}{\eta}}} & &
\end{eqnarray*}
and as a consequence
\begin{eqnarray*}
 \int_{0}^{1} \sum_{i=1}^{+\infty} \P( \sup_{w \in B(u,\varep)} [|f(v,w,V_1)\log|f(v,w,V_1)||^{\alpha(w)}]  > \frac{Kx^{\frac{d}{\eta}} }{d^d}i) dx \\
& \leq & \int_{0}^{1} \frac{K'}{x^{\frac{d}{\eta}}} dx \\
& < & +\infty .\\
\end{eqnarray*}

\smallskip
\smallskip

We have thus verified all the conditions in the three series theorem, and
shown that the series $Z^1, Z^2$ and $Z^3$ are almost surely convergent.

\medskip

{\bf Second step: almost sure convergence of $\sum_{i=1}^{\infty}\limits Y_i^j, j=1,2,3$.}

\medskip

To prove that the series $\sum_{i=1}^{\infty}\limits Y_i^j, j=1,2,3$ converge almost surely, we will first show that it is enough to prove that $\sum_{i=1}^{\infty}\limits Y_i^j \mathbf{1}_{\{\frac{1}{2} \leq \frac{\Gamma_i}{i} \leq 2\} \cap \{|Y_i^j| \leq 1\}}$ converges almost surely
for $j=1,2,3$. Indeed, we prove now that
$\sum_{i=1}^{\infty}\P \left( \overline{\{ \frac{1}{2} \leq \frac{\Gamma_i}{i} \leq 2\}}  \cup \{|Y_i^j| >1 \} \right) < \infty$, where $\overline T$ denotes
the complementary set of the set $T$.

 \begin{eqnarray*}
  \P \left( \overline{\{ \frac{1}{2} \leq \frac{\Gamma_i}{i} \leq 2\}}  \cup \{|Y_i^j| >1 \} \right) & = & \P\left( \overline{\{ \frac{1}{2} \leq \frac{\Gamma_i}{i} \leq 2\}} \cup  \left[ \{ |Y_i^j| >1\} \cap \{ \frac{1}{2} \leq \frac{\Gamma_i}{i} \leq 2\} \right] \right)\\
  & \leq & \P(\{  \Gamma_i < \frac{i}{2}\}) + \P(\{ \Gamma_i > 2i\}) + \P \left( \{ |Y_i^j| >1\} \cap \{ \frac{1}{2} \leq \frac{\Gamma_i}{i} \leq 2\} \right).\\
 \end{eqnarray*}
$\Gamma_i$, as a sum of independent and identically distributed exponential random variables with mean $1$, satisfy a Large Deviation Principle with rate function $\Lambda^*(x)=x-1-\log(x)$ for $x>0$ and infinity for $x \leq 0$ (see for instance \cite{DZ} p.35), thus $\sum_{i \geq 1}\limits \P(\{  \Gamma_i < \frac{i}{2}\}) < + \infty$ and $\sum_{i \geq 1}\limits \P(\{ \Gamma_i > 2i\}) < + \infty$.

Consider now $\sum_{i \geq 1}\limits \P \left(\{ |Y_i^j| >1\} \cap \{ \frac{1}{2} \leq \frac{\Gamma_i}{i} \leq 2\}\right),$ for $j=1,2,3$.

\medskip

\underline {Case $j=1$ or $j=2$} :

\smallskip

\begin{eqnarray*}
 \P \left(\{ |Y_i^j| >1\} \cap \{ \frac{1}{2} \leq \frac{\Gamma_i}{i} \leq 2\}\right) & =& \P \left( \{ |Z_i^j||(\frac{\Gamma_i}{i})^{-1/ \alpha(x_i)}-1|>1\} \cap \{\frac{1}{2} \leq \frac{\Gamma_i}{i} \leq 2 \}  \right) \\
 & \leq & \P \left( \{ (2^{1/ \alpha(x_i)}-1)|Z_i^j|>1\} \cap \{\frac{1}{2} \leq \frac{\Gamma_i}{i} \leq 2 \} \right)\\
 & \leq & \P \left( |Z_i^j| > \frac{1}{2^{1/d}-1}\right) \\
\end{eqnarray*}
thus $\sum_{i \geq 1} \P \left(\{ |Y_i^j| >1\} \cap \{ \frac{1}{2} \leq \frac{\Gamma_i}{i} \leq 2\}\right) < + \infty.$

\medskip

\underline {Case $j=3$}:

\smallskip
 For $i >1:$
\begin{eqnarray*}
 \P \left(\{ |Y_i^3| >1\} \cap \{ \frac{1}{2} \leq \frac{\Gamma_i}{i} \leq 2\}\right) & =& \P \left( \{ |Z_i^3||\frac{\log \Gamma_i}{\log i}(\frac{\Gamma_i}{i})^{-1/ \alpha(x_i)}-1|>1\} \cap \{\frac{1}{2} \leq \frac{\Gamma_i}{i} \leq 2 \}  \right) .\\
\end{eqnarray*}

$(\frac{\log \Gamma_i}{\log i}(\frac{\Gamma_i}{i})^{-1/ \alpha(x_i)}-1)_{(i>1)}$ is bounded on $\{\frac{1}{2} \leq \frac{\Gamma_i}{i} \leq 2  \}$, thus there exists $K>0$ such that
\begin{displaymath}
 \P \left(\{ |Y_i^3| >1\} \cap \{ \frac{1}{2} \leq \frac{\Gamma_i}{i} \leq 2\}\right) \leq \P\left( |Z_i^3| > \frac{1}{K} \right)
\end{displaymath}
which entails that $\sum_{i \geq 1}\limits \P \left(\{ |Y_i^3| >1\} \cap \{ \frac{1}{2} \leq \frac{\Gamma_i}{i} \leq 2\}\right) < + \infty.$

\medskip

We are thus left with proving that $ \sum_{i=1}^{\infty}\limits Y_i^j \mathbf{1}_{\{\frac{1}{2} \leq \frac{\Gamma_i}{i} \leq 2\} \cap \{|Y_i^j| \leq 1\}}$ converges almost surely  for $j=1,2,3$.

In that view, we shall apply the following well-known lemma:

\begin{lem}
 Let $\{ X_k, k\geq 1\}$ be a sequence of random variables such that $\sum_{n=1}^{+\infty}\limits \E |X_n| < + \infty$, then $\sum_{n=1}^{+\infty}\limits X_n$ converges almost surely.
\end{lem}

%

Let us show that $\sum_{i=1}^{\infty}\limits \E \left[|Y_i^j| \mathbf{1}_{\{\frac{1}{2} \leq \frac{\Gamma_i}{i} \leq 2\} \cap \{|Y_i^j| \leq 1\}} \right] < +\infty$.

\begin{eqnarray*}
 \E \left[ |Y_i^j| \mathbf{1}_{\{\frac{1}{2} \leq \frac{\Gamma_i}{i} \leq 2\} \cap \{|Y_i^j| \leq 1\}} \right] &= & \int_{0}^{\infty} \P \left( \{ 1 \geq |Y_i^j| > x \} \cap \{ \frac{1}{2} \leq \frac{\Gamma_i}{i} \leq 2\}\right) dx \\
 & \leq &\int_{0}^{1} \P \left(\{|Y_i^j| > x \} \cap \{ \frac{1}{2} \leq \frac{\Gamma_i}{i} \leq 2\} \right) dx .\\
\end{eqnarray*}

Let $B_i =\{ \frac{1}{2} \leq \frac{\Gamma_i}{i} \leq 2\}$.
\smallskip

\underline {Case $j=1$}:

Using the finite-increments formula applied to the function $x \mapsto x^{-\frac{1}{\alpha(x_i)}}$ on $[\frac{1}{2},2]$, one easily shows that

 \begin{eqnarray*}
  \P \left(\{|Y_i^1| > x \} \cap B_i \right) & \leq & \P \left( \{ |Z_i^1| |\frac{\Gamma_i}{i}-1|>x \frac{c}{2^{1+1/c}} \} \cap B_i \right) \\
 & = & \P \left( |a'(x_i) i^{-1/ \alpha(x_i)} f(v,x_i,V_i)| |\frac{\Gamma_i}{i}-1|>x \frac{c}{2^{1+1/c}} \} \cap B_i \right) \\
 & \leq & \P \left( \{ |f(v,x_i,V_i)|^{\alpha(x_i)} |\frac{\Gamma_i}{i}-1|^{\alpha(x_i)} > K_c i x^{\alpha(x_i)} \} \cap B_i \right)\\
 \end{eqnarray*}
 where $K_c:= \inf_{w \in B(u,\varep)}
\left[\left(\frac{c}{2^{1+1/c}|a'(w)|} \right)^{\alpha(w)}\right]$ is strictly positive by the assumptions on $a'$ and $\alpha$.  Thus
 \begin{eqnarray*}
 \P \left(\{|Y_i^1| > x \} \cap B_i \right) & \leq & \P \left( \sup_{w \in B(u,\varep)} |f(v,w,V_1)|^{\alpha(w)} |\frac{\Gamma_i}{i}-1|^{c} > K_c i x^{d} \right).\\
 \end{eqnarray*}

 \underline{{\it Case $d\geq 1$}}:

 \smallskip

 Fix $\eta \in (d,1+\frac{c}{2})$ (since $\alpha$ is continuous and
 $d<2$, by decreasing if necessary $\varepsilon$, one may ensure that
 $d < 1+ c/2$). By Markov and H\"{o}lder inequalities,

 \begin{eqnarray*}
  \P \left(\{|Y_i^1| > x \} \cap B_i \right) & \leq & \P \left( \left[\sup_{w \in B(u,\varep)} |f(v,w,V_1)|^{\alpha(w)} \right]^{1/ \eta} |\frac{\Gamma_i}{i}-1|^{c/ \eta} > K_c^{1/ \eta} i^{1/ \eta} x^{d/\eta} \right)\\
  & \leq & \frac{1}{(K_cix^d)^{1/ \eta}} \left[ \E |\frac{\Gamma_i}{i}-1|^2 \right]^{c/2\eta} \left( \sup_{v \in B(u,\varep)} \E (\sup_{w \in B(u,\varep)} |f(v,w,V_1)|^{\alpha(w)})\right)^{1/ \eta}\\
  & \leq & \frac{K}{x^{d/\eta}}\frac{1}{i^{1/\eta+c/2\eta}}\\
 \end{eqnarray*}
 where we have used that the variance of $\Gamma_i$ is equal to $i$, and
 K does not depend on $v$ thanks to assumption (C2).
Thus $\E \left[ |Y_i^1| \mathbf{1}_{\{\frac{1}{2} \leq \frac{\Gamma_i}{i} \leq 2\} \cap \{|Y_i^1| \leq 1\}} \right] \leq \frac{K}{i^{1/\eta+c/2\eta}}$ where $\frac{1}{\eta} + \frac{c}{2\eta} > 1$.

\smallskip

\underline{{\it Case $d < 1$}}:

 \smallskip

 \begin{eqnarray*}
   \P \left(\{|Y_i^1| > x \} \cap B_i \right) & \leq & \frac{1}{x^d K_c i} \E (\sup_{w \in B(u,\varep)} |f(v,w,V_1)|^{\alpha(w)}) \E |\frac{\Gamma_i}{i}-1|^c\\
   & \leq & K \frac{1}{x^d} \frac{1}{i} (\E |\frac{\Gamma_i}{i}-1|^2)^{c/2}\\
   & \leq & K \frac{1}{i^{1+c/2}} \frac{1}{x^d},\\
 \end{eqnarray*}
thus $\E \left[ |Y_i^1| \mathbf{1}_{\{\frac{1}{2} \leq \frac{\Gamma_i}{i} \leq 2\} \cap \{|Y_i^1| \leq 1\}} \right] \leq \frac{K}{i^{1 + c/2}}$ with $1 + \frac{c}{2} > 1$.

\bigskip

The case of $\sum_{i \geq 1}\limits \E \left[|Y_i^2| \mathbf{1}_{\{B_i \cap \{|Y_i^2| \leq 1\}} \right]$ is treated similarly, since the conditions required on $(a',f)$ in the proof above are also satisfied by $(a,f'_u)$.

\bigskip

\underline {Case $j=3$}:

We now consider $\sum_{i \geq 1}\limits \E \left[|Y_i^3| \mathbf{1}_{\{B_i \cap \{|Y_i^3| \leq 1\}} \right]$.

\smallskip

Again by the finite-increments formula, there exists a constant $K_{c,d}$, which depends on $c$ and $d$, such that, for $i>1$,
\begin{displaymath}
 \left| \frac{\log \Gamma_i}{\log i}(\frac{\Gamma_i}{i})^{-1/ \alpha(x_i)}-1\right| \mathbf{1}_{B_i} \leq K_{c,d} |\frac{\Gamma_i}{i}-1| \mathbf{1}_{B_i}.
\end{displaymath}

Then,

\begin{eqnarray*}
  \P \left(\{|Y_i^3| > x \} \cap B_i \right) & \leq & \P \left( \{ |Z_i^3| K_{c,d} |\frac{\Gamma_i}{i}-1| > x \} \cap B_i\right) \\
  & \leq & \P \left(  \sup_{w \in B(u,\varep)} |f(v,w,V_1)|^{\alpha(w)} |\frac{\Gamma_i}{i}-1|^c > K \frac{i}{(\log i)^d} x^d \right).\\
\end{eqnarray*}

\underline{{\it Case $d\geq 1$}} :

 \smallskip

 Fix $\eta \in (d,1+\frac{c}{2})$.

 \begin{eqnarray*}
  \P \left(\{|Y_i^3| > x \} \cap B_i \right) & \leq & \frac{K}{x^{d/\eta}}\frac{(\log i)^{d/\eta}}{i^{1/\eta+c/2\eta}}.\\
 \end{eqnarray*}

\underline{{\it Case $d < 1$}} :

 \smallskip

 \begin{eqnarray*}
   \P \left(\{|Y_i^3| > x \} \cap B_i \right) & \leq & \frac{K}{x^d} \frac{(\log i)^d}{i^{1+c/2}}.\\
 \end{eqnarray*}

\bigskip

As a conclusion, for $j=1,2,3$, $\sum_{i=1}^{+\infty}\limits Y_i^j$ converges almost surely.

\medskip

We now move to the last two steps of the proof: to verify $h$-localisability, we need to check that for some $\eta$ such that $h<\eta<1$,
$\P(|Z^j|\geq |v-u|^{\eta-1})$ and $\P(|Y^j|\geq |v-u|^{\eta-1})$ tend to 0 when $v$ tends to $u$, for $j=1,2,3$.

\medskip

{\bf Third step: verification of (\ref{cond}) for $Z^j, j=1,2,3$.}

\medskip
We need to estimate $\P(|Z^j|\geq |v-u|^{\eta-1})$.

Let $a \in (0, \frac{2(1-\eta)}{2-c})$.

\begin{eqnarray*}
 \P\left( |\sum_{i=1}^{\infty} Z_i^j | > |v-u|^{\eta-1} \right) & \leq & \P \left( |\sum_{i=1}^{\infty} Z_i^j \mathbf{1}_{|Z_i^j| \leq |v-u|^{-a}} | > \frac{|v-u|^{\eta-1}}{2} \right) \\
&  + & \P \left( |\sum_{i=1}^{\infty} Z_i^j \mathbf{1}_{|Z_i^j| > |v-u|^{-a}} | > \frac{|v-u|^{\eta-1}}{2}\right).
\end{eqnarray*}

Since $\gamma_i$ is independent from $\gamma_k$ for $i \neq k$ and $|Z_i^j |$ is independent of
$\gamma_i$, Markov inequality yields

\begin{displaymath}
 \P \left( |\sum_{i=1}^{\infty} Z_i^j \mathbf{1}_{|Z_i^j| \leq |v-u|^{-a}} | > \frac{|v-u|^{\eta-1}}{2} \right) \leq \frac{4}{|v-u|^{2(\eta-1)}} \sum_{i=1}^{\infty} \E \left[ |Z_i^j|^2 \mathbf{1}_{|Z_i^j| \leq |v-u|^{-a}}\right].
\end{displaymath}

Let $\gamma \in (d,2)$. For any $M>1$, using the same computations as in the first step, third series (page \pageref{page10}), we get:

\begin{displaymath}
 \E \left[ |Z_i^j|^2 \mathbf{1}_{|Z_i^j| \leq M}\right] = M^2\E \left[ \frac{|Z_i^j|^2}{M^2} \mathbf{1}_{|Z_i^j| \leq M}\right] \leq M^2\int_{0}^{1} \P(|Z_i^j| > M x^{1/ \gamma} )dx.
\end{displaymath}

For $j=1$,

\begin{displaymath}
 \int_{0}^{1} \P(|Z_i^1| > M x^{1/ \gamma} )dx \leq \int_{0}^{1} \P(\sup_{w \in B(u,\varep)}| f(v,w,V_1)|^{\alpha(w)} > M^c i K x^{d/\gamma}) dx,
\end{displaymath}
thus
\begin{displaymath}
 \sum_{i=1}^\infty \E \left[ |Z_i^1|^2 \mathbf{1}_{|Z_i^1| \leq M}\right] \leq M^2\sum_{i=1}^\infty \int_{0}^{1} \P(\sup_{w \in B(u,\varep)}| f(v,w,V_1)|^{\alpha(w)} > M^c i K x^{d/\gamma}) dx \leq K M^{2-c}.
\end{displaymath}

The same conclusion holds for $j=2$:
\begin{displaymath}
\sum_{i=1}^\infty \E \left[ |Z_i^2|^2 \mathbf{1}_{|Z_i^2| \leq M}\right] \leq K M^{2-c}.
\end{displaymath}

For $j=3$, choose $\mu \in (0,1-\frac{d}{\gamma})$, and compute as on page \pageref{page11}:

\begin{eqnarray*}
 \int_{0}^{1} \P(|Z_i^3| > M x^{1/ \gamma}) dx & \leq & \int_{0}^{1} \P ( A > |f(v,w_i,V_i)|^{\alpha(w_i)} > K M^c x^{d/ \gamma} \frac{i}{\log(i)^d})dx\\
 & & + \int_{0}^{1} \P( \sup_{w \in B(u,\varep)} [| f(v,w,V_1)\log f(v,w,V_1)| ]^{\alpha(w)} > K \frac{M^c x^{d/\eta}}{d^d} i ) dx.
\end{eqnarray*}

For $i$ large enough, $i > i^*$ where $i^*$ depends only on $d$ and $\mu$,

\begin{displaymath}
 \P ( A > |f(v,w_i,V_i)|^{\alpha(w_i)} > K M^c x^{d/ \gamma} \frac{i}{\log(i)^d}) \leq \P ( A > |f(v,w_i,V_i)|^{\alpha(w_i)} > K M^c x^{d/ \gamma} i^{1-\mu})
\end{displaymath}
thus

\begin{eqnarray*}
 \sum_{i=1}^{\infty} \P \left( A > |f(v,w_i,V_i)|^{\alpha(w_i)} > K M^c x^{d/ \gamma} \frac{i}{\log(i)^d} \right) & \leq & \\
   \left( \sum_{i=1}^{i^*} \frac{\log(i)^d}{i}\right) \frac{x^{-d / \gamma}}{KM^c} \sup_{v \in B(u,\varep)} \E (\sup_{w \in B(u,\varep)} | f(v,w,V_1)|^{\alpha(w)}) &  &\\
+ \sum_{i=i*}^{\infty} \P ( A > |f(v,w_i,V_i)|^{\alpha(w_i)} > K M^c x^{d/ \gamma} i^{1-\mu}).\\
\end{eqnarray*}

As a consequence, for $M > 1$:
\begin{eqnarray*}
 \int_0^1 \sum_{i=1}^{\infty} \P \left( A > |f(v,w_i,V_i)|^{\alpha(w_i)} > K M^c x^{d/ \gamma} \frac{i}{\log(i)^d} \right) & \leq & K M^{-c} + K M^{-\frac{c}{1-\mu}}\\
& \leq & K M^{-c}.\\
\end{eqnarray*}
 Since
\begin{displaymath}
 \sum_{i=1}^{\infty} \int_{0}^{1} \P( \sup_{w \in B(u,\varep)} [| f(v,w,V_1)\log |f(v,w,V_1)|| ]^{\alpha(w)} > K \frac{M^c x^{d/\eta}}{d^d} i ) dx \leq K M^{-c},
\end{displaymath}
we get
\begin{displaymath}
 \sum_{i=1}^\infty \E \left[ |Z_i^3|^2 \mathbf{1}_{|Z_i^3| \leq M}\right] \leq K M^{2-c}.
\end{displaymath}

Let $M = |v-u|^{-a}$. Using previously obtained inequalities, we get, for $j=1,2,3$:
\begin{displaymath}
 \P \left( |\sum_{i=1}^{\infty} Z_i^j \mathbf{1}_{|Z_i^j| \leq |v-u|^{-a}} | > \frac{|v-u|^{\eta-1}}{2} \right) \leq K |v-u|^{2(1-\eta)-a(2-c)}
\end{displaymath}
and
\begin{displaymath}
 \lim_{v \rightarrow u} \P \left( |\sum_{i=1}^{\infty} Z_i^j \mathbf{1}_{|Z_i^j| \leq |v-u|^{-a}} | > \frac{|v-u|^{\eta-1}}{2} \right) =0.
\end{displaymath}
We consider now the second term $\P \left( |\sum_{i=1}^{\infty} Z_i^j \mathbf{1}_{|Z_i^j| > |v-u|^{-a}} | > \frac{|v-u|^{\eta-1}}{2} \right) $.

Let $i^{*}=\inf\{n\geq1: i\geq n,|Z_i^j|\leq |v-u|^{-a} \}$. Since $\sum_{i \geq 1}\limits \P (|Z_i^j| > |v-u|^{-a} ) < +\infty$,
the Borel-Cantelli lemma yields $\P (i^* = + \infty) = 0$. As a consequence,

\begin{eqnarray*}
 \P \left( |\sum_{i=1}^{\infty} Z_i^j \mathbf{1}_{|Z_i^j| > |v-u|^{-a}} | > \frac{|v-u|^{\eta-1}}{2} \right) & = & \sum_{n=1}^{\infty} \P \left( \{|\sum_{i=1}^{\infty} Z_i^j \mathbf{1}_{|Z_i^j| > |v-u|^{-a}} | > \frac{|v-u|^{\eta-1}}{2} \} \cap \{ i^*=n \}\right) \\
& = & \sum_{n=2}^{\infty} \P \left( \{|\sum_{i=1}^{\infty} Z_i^j \mathbf{1}_{|Z_i^j| > |v-u|^{-a}} | > \frac{|v-u|^{\eta-1}}{2} \} \cap \{ i^*=n \}\right) \\
 & \leq &  \sum_{n=2}^{\infty} \P (i^*=n).\\
\end{eqnarray*}

For $n \geq 2$, $\P (i^*=n) \leq \P( |Z_{n-1}^j| > |v-u|^{-a})$.

For $j=1$, $\P (i^*=n) \leq \P( \sup_{w \in B(u,\varep)} |f(v,w,V_1)|^{\alpha(w)} > |v-u|^{-ac} K (n-1))$,
and thus

\begin{eqnarray*}
  \sum_{n=2}^{\infty} \P (i^*=n) & \leq  &K |v-u|^{ac} \E ( \sup_{w \in B(u,\varep)} |f(v,w,V_1)|^{\alpha(w)})\\
& \leq & K |v-u|^{ac} \sup_{t \in B(u,\varep)} \E ( \sup_{w \in B(u,\varep)} |f(t,w,V_1)|^{\alpha(w)}).\\
\end{eqnarray*}

For $j=2$,
\begin{displaymath}
\sum_{n=2}^{\infty} \P (i^*=n) \leq  K |v-u|^{ac} \sup_{t \in B(u,\varep)} \E ( \sup_{w \in B(u,\varep)} |f'_u(t,w,V_1)|^{\alpha(w)}),
\end{displaymath}
and for $j=3$,
\begin{eqnarray*}
 \sum_{n=2}^{\infty} \P (i^*=n) & \leq &  K |v-u|^{ac} \sup_{t \in B(u,\varep)} \E ( \sup_{w \in B(u,\varep)}|f(t,w,V_1)\log |f(t,w,V_1)||^{\alpha(w)} )\\
& + & \sum_{n\geq 2} \P( A > |f(v,w_i,V_i)|^{\alpha(w_i)} > |v-u|^{-a\alpha(w_i)} \frac{K i }{\log(i)^d}  ).
\end{eqnarray*}

We have shown previously that the second term in the sum on the right hand side of the above inequality
is bounded from above by $K |v-u|^{ac}$.
Finally,
\begin{displaymath}
  \lim_{v \rightarrow u} \P \left( |\sum_{i=1}^{\infty} Z_i^j \mathbf{1}_{|Z_i^j| > |v-u|^{-a}} | > \frac{|v-u|^{\eta}}{2} \right) =0.
\end{displaymath}

\medskip

{\bf Fourth step: verification of (\ref{cond}) for $Y^j, j=1,2,3$.}

\medskip

We consider now $\P(|Y^j|\geq |v-u|^{\eta-1})$.

Let $i^{*}=\inf\{n\geq1: i\geq n,|Y_i^j|\leq 1 \textrm{ and } \frac{1}{2} \leq \frac{\Gamma_i}{i} \leq 2\}$. Since $\sum_{i \geq 1}\limits \P ( \{|Y_i^j| > 1 \} \cup \{ \Gamma_i < \frac{i}{2}\} \cup \{ \Gamma_i > 2i \} ) < +\infty$,
the Borel-Cantelli lemma yields $\P (i^* = + \infty) = 0$. As a consequence,
\begin{eqnarray*}
 \P(|Y^j| \geq |v-u|^{\eta - 1}) & = & \sum_{n \geq 1} \P \left(\{ |\sum_{i=1}^{\infty} Y_i^j| \geq |v-u|^{\eta-1}\} \cap \{ i^{*} =n\} \right).\\
 \end{eqnarray*}
Let $b_n(v) = \P \left(\{ |\sum_{i=1}^{\infty}\limits Y_i^j| \geq |v-u|^{\eta-1}\} \cap \{ i^{*} =n\} \right)$.
Our strategy is the following: we show that, for each fixed $n$, $b_n(v)$ tend to 0 when $v$ tends to $u$. Then
we prove that there exists a summable sequence $(c_n)_n$ such that, for all $n$ and all $v$, $b_n(v) \leq c_n$. We conclude using the dominated convergence theorem that $\sum_{n \geq 1} b_n(v)$ tends to 0 when $v$ tends to $u$.

For all $n \geq 1$,
\begin{eqnarray*}
  b_n(v) & \leq & \P (\{ |\sum_{i=1}^{n-1} Y_i^j| \geq \frac{|v-u|^{\eta-1}}{2} \}) +  \P (\{ |\sum_{i=n}^{\infty} Y_i^j\mathbf{1}_{\{|Y_i^j|\leq 1 \} \cap \{ \frac{1}{2} \leq \frac{\Gamma_i}{i} \leq 2\}}| \geq \frac{|v-u|^{\eta-1}}{2}\}).\\
 \end{eqnarray*}
For $n \geq 2$, consider $\P (\{ |\sum_{i=1}^{n-1}\limits Y_i^j| \geq \frac{|v-u|^{\eta-1}}{2} \})$.

\begin{eqnarray*}
\P (\{ |\sum_{i=1}^{n-1} Y_i^j| \geq \frac{|v-u|^{\eta-1}}{2} \}) &\leq & \sum_{i=1}^{n-1} \P (| Y_i^j | \geq \frac{|v-u|^{\eta-1}}{2(n-1)}).\\
\end{eqnarray*}

Let $p \in (0,\frac{c}{d}).$ With $K$ a positive constant that may change from line to line and depend on $n$ but not on $v$, we have, for $j=1$:

\begin{eqnarray*}
\P \left(| Y_i^1 | \geq \frac{|v-u|^{\eta-1}}{2(n-1)}\right) & \leq & \P \left( \sup_{w \in B(u,\varep)} |f(v,w,V_i)|^{\alpha(w)} |(\frac{\Gamma_i}{i})^{-1/ \alpha(x_i)}-1|^{\alpha(x_i)} \geq \frac{i|v-u |^{\alpha(x_i)(\eta-1)}}{(2(n-1)a'(x_i))^{\alpha(x_i)}}  \right) \\
& \leq & \P \left( \sup_{w \in B(u,\varep)} |f(v,w,V_i)|^{\alpha(w)}  |(\frac{\Gamma_i}{i})^{-1/ \alpha(x_i)}-1|^{\alpha(x_i)}\geq K |v-u|^{c(\eta-1)} \right)\\
& \leq & \P \left( (\sup_{w \in B(u,\varep)} |f(v,w,V_i)|^{\alpha(w)})^p  |(\frac{\Gamma_i}{i})^{-1/ \alpha(x_i)}-1|^{p\alpha(x_i)}\geq K |v-u|^{pc(\eta-1)} \right).\\
\end{eqnarray*}

Using the independence of $V_i$ and $\Gamma_i$ and Markov inequality,

\begin{eqnarray*}
\P \left(| Y_i^1 | \geq \frac{|v-u|^{\eta-1}}{2(n-1)}\right) & \leq & K|v-u|^{pc(1-\eta)} \E\left((\sup_{w \in B(u,\varep)} |f(v,w,V_i)|^{\alpha(w)})^p\right) \E(|(\frac{\Gamma_i}{i})^{-1/ \alpha(x_i)}-1|^{p\alpha(x_i)})\\
& \leq & K |v-u|^{pc(1-\eta)} \E (|(\frac{\Gamma_i}{i})^{-1/ \alpha(x_i)}-1|^{p\alpha(x_i)}).
\end{eqnarray*}
It remains to check that the expectation in the right hand side of the above inequality is finite:
\begin{eqnarray*}
\E (|(\frac{\Gamma_i}{i})^{-1/ \alpha(x_i)}-1|^{p\alpha(x_i)}) & = & \E (|(\frac{\Gamma_i}{i})^{-1/ \alpha(x_i)}-1|^{p\alpha(x_i)}\mathbf{1}_{\Gamma_i >i}) + \E (|(\frac{\Gamma_i}{i})^{-1/ \alpha(x_i)}-1|^{p\alpha(x_i)}\mathbf{1}_{\Gamma_i <i})\\
& \leq & 1+\E (|(\frac{\Gamma_i}{i})^{-1/ \alpha(x_i)}-1|^{p\alpha(x_i)}\mathbf{1}_{\Gamma_i <i})\\
& \leq & 1 + \E (|(\frac{\Gamma_i}{i})^{-1/ c}-1|^{p\alpha(x_i)}\mathbf{1}_{\Gamma_i <i})\\
& \leq & 1+ \E (|(\frac{\Gamma_i}{i})^{-1/ c}-1|^{pc})+ \E (|(\frac{\Gamma_i}{i})^{-1/ c}-1|^{pd}).\\
\end{eqnarray*}
Since $p < \frac{c}{d}$, $\E (|(\frac{\Gamma_i}{i})^{-1/ c}-1|^{pc}) < + \infty$ and $\E(|(\frac{\Gamma_i}{i})^{-1/ c}-1|^{pd}) < + \infty$ (this is easily verified by computing these expectations using the density of $\Gamma_i$). Thus we have $\E(|(\frac{\Gamma_i}{i})^{-1/ \alpha(x_i)}-1|^{p\alpha(x_i)}) < +\infty$, and
\begin{displaymath}
\lim_{v \rightarrow u} \P \left(| Y_i^1 | \geq \frac{|v-u|^{\eta-1}}{2(n-1)}\right)=0.
\end{displaymath}
Since the conditions required on $(a',f)$ are also satisfied by $(a,f'_u)$,
\begin{displaymath}
\lim_{v \rightarrow u} \P \left(| Y_i^2 | \geq \frac{|v-u|^{\eta-1}}{2(n-1)}\right)=0.
\end{displaymath}

We consider now the case $j=3$. When $i=1$:
\begin{eqnarray*}
\P \left( |Y_1^3| \geq \frac{|v-u|^{\eta-1}}{2(n-1)}\right) & =&  \P \left( |\log(\Gamma_1) \Gamma_1^{-1/ \alpha(x_1)} f(v,x_1,V_1) | \geq \frac{\alpha(x_1)^2}{2 |a(x_1) \alpha'(x_1)|(n-1)} |v-u|^{\eta-1} \right)\\
& \leq & K |v-u|^{pc(1-\eta)} \E \left( ( \frac{(\log(\Gamma_1))^{\alpha(x_1)}}{\Gamma_1})^p\right),
\end{eqnarray*}
where again $K$ depends on $n$ but not on $v$.
Since $p<1$ and $\alpha$ is bounded, $\E \left( ( \frac{(\log(\Gamma_1))^{\alpha(x_1)}}{\Gamma_1})^p\right) < + \infty$,  and
\begin{displaymath}
\lim_{v \rightarrow u}\limits \P \left(| Y_1^3 | \geq \frac{|v-u|^{\eta-1}}{2(n-1)}\right)=0.
\end{displaymath}

For $i \geq 2$,
\begin{eqnarray*}
\P \left(| Y_i^3 | \geq \frac{|v-u|^{\eta-1}}{2(n-1)}\right) & = & \P \left(\left|\left(\frac{\log(\Gamma_i)}{\log(i)}(\frac
{\Gamma_i}{i})^{-1/\alpha(x_i)}-1\right) f(v,x_i,V_i)\right| \geq \frac{\alpha(x_i)^2 i^{1/\alpha(x_i)}|v-u|^{\eta-1}}{\log(i) 2 |a(x_i) \alpha'(x_i)| (n-1)} \right)\\
& \leq & K\left(\frac{\log(i)^d}{i}\right)^p |v-u|^{pc(1-\eta)} \E \left( | \frac{\log(\Gamma_i)}{\log(i)}(\frac
{\Gamma_i}{i})^{-1/\alpha(x_i)}-1|^{\alpha(x_i)p}\right).
\end{eqnarray*}
\begin{eqnarray*}
\E \left( | \frac{\log(\Gamma_i)}{\log(i)}(\frac{\Gamma_i}{i})^{-1/\alpha(x_i)}-1|^{\alpha(x_i)p}\right) & \leq & \E\left( | \frac{\log(\Gamma_i)}{\log(i)}(\frac{\Gamma_i}{i})^{-1/c}-1|^{\alpha(x_i)p} +  | \frac{\log(\Gamma_i)}{\log(i)}(\frac{\Gamma_i}{i})^{-1/d}-1|^{\alpha(x_i)p}\right).
\end{eqnarray*}
\begin{eqnarray*}
\E \left( | \frac{\log(\Gamma_i)}{\log(i)}(\frac{\Gamma_i}{i})^{-1/c}-1|^{\alpha(x_i)p}\right) & \leq & \E \left( | \frac{\log(\Gamma_i)}{\log(i)}(\frac{\Gamma_i}{i})^{-1/c}-1|^{cp}\right)+\E \left( | \frac{\log(\Gamma_i)}{\log(i)}(\frac{\Gamma_i}{i})^{-1/c}-1|^{dp}\right)\\
\end{eqnarray*}
and
\begin{eqnarray*}
\E \left( | \frac{\log(\Gamma_i)}{\log(i)}(\frac{\Gamma_i}{i})^{-1/d}-1|^{\alpha(x_i)p}\right) & \leq & \E \left( | \frac{\log(\Gamma_i)}{\log(i)}(\frac{\Gamma_i}{i})^{-1/d}-1|^{cp}\right)+\E \left( | \frac{\log(\Gamma_i)}{\log(i)}(\frac{\Gamma_i}{i})^{-1/d}-1|^{dp}\right).\\
\end{eqnarray*}

Since $p \in (0,\frac{c}{d})$, the four terms in the right hand sides of the two last inequalities are finite (use again the density of $\Gamma_i$) and thus $\E \left(|\frac{\log(\Gamma_i)}{\log(i)}(\frac{\Gamma_i}{i})^{-1/\alpha(x_i)}-1|^{\alpha(x_i)p}\right) < + \infty$. As a consequence,
\begin{displaymath}
\lim_{v \rightarrow u} \P \left(| Y_i^3 | \geq \frac{|v-u|^{\eta-1}}{2(n-1)}\right)=0.
\end{displaymath}
Finally, we have, for $j \in \{ 1,2,3\}$,
\begin{displaymath}
\lim_{v \rightarrow u}\P (\{ |\sum_{i=1}^{n-1} Y_i^j| \geq \frac{|v-u|^{\eta-1}}{2} \})=0.
\end{displaymath}

Let us now consider, for $n \geq 1$, $\P \left(\{ |\sum_{i=n}^{\infty}\limits Y_i^j\mathbf{1}_{\{|Y_i^j|\leq 1\} \cap \{ \frac{1}{2} \leq \frac{\Gamma_i}{i} \leq 2\}}| \geq \frac{|v-u|^{\eta-1}}{2}\}\right)$:
\begin{eqnarray*}
\P \left(\{ |\sum_{i=n}^{\infty} Y_i^j\mathbf{1}_{\{|Y_i^j|\leq 1\} \cap \{ \frac{1}{2} \leq \frac{\Gamma_i}{i} \leq 2\}}| \geq \frac{|v-u|^{\eta-1}}{2}\}\right) & \leq &  2 |v-u|^{1-\eta} \E \left[ |\sum_{i=n}^{\infty} Y_i^j\mathbf{1}_{\{|Y_i^j|\leq 1\} \cap \{ \frac{1}{2} \leq \frac{\Gamma_i}{i} \leq 2\}}|\right]\\
& \leq &  2 |v-u|^{1-\eta} \sum_{i=1}^{\infty} \E|Y_i^j|\mathbf{1}_{\{|Y_i^j|\leq 1\} \cap \{ \frac{1}{2} \leq \frac{\Gamma_i}{i} \leq 2\}} \\
& \leq & K |v-u|^{1-\eta}\\
\end{eqnarray*}
(recall that the constants $K$ used in bounding the series $\E(|Y_i^j|\mathbf{1}_{\{|Y_i^j|\leq 1\} \cap \{ \frac{1}{2} \leq \frac{\Gamma_i}{i} \leq 2\}})$ do not depend on $v$). Thus $b_n(v) \to 0$
when $v \to u$ for each $n$.

In view of using the dominated convergence theorem, we compute (recall that $B_i = \{ \frac{1}{2} \leq \frac{\Gamma_i}{i} \leq 2 \}$):
\begin{eqnarray*}
b_n(v) & \leq & \P(\{i^*=n\})\\
& \leq & \P( \{ | Y_{n-1}^j| > 1\} \cup \overline{B_{n-1}})\\
& \leq & \P( \{ | Y_{n-1}^j| > 1\} \cap B_{n-1}) + \P( \frac{\Gamma_{n-1}}{n-1} < \frac{1}{2} ) + \P ( \frac{\Gamma_{n-1}}{n-1} > 2) .\\
\end{eqnarray*}
For $j=1$ and $d\geq1$,

\begin{displaymath}
 \P( \{ | Y_{n-1}^1| > 1\} \cap B_{n-1}) \leq \frac{K}{(n-1)^{1/ \eta +c/2\eta}} ( \sup_{t \in B(u,\varep)} \E ( \sup_{w \in B(u,\varep)} |f(t,w,V_1)|^{\alpha(w)}) )^{1/ \eta}
\end{displaymath}
and if $d < 1$,
\begin{displaymath}
 \P( \{ | Y_{n-1}^1| > 1\} \cap B_{n-1}) \leq \frac{K}{(n-1)^{1 +c/2}} ( \sup_{t \in B(u,\varep)} \E ( \sup_{w \in B(u,\varep)} |f(t,w,V_1)|^{\alpha(w)}) ).
\end{displaymath}
The same conclusion holds for $j=2$, while, for $j=3$,
\begin{displaymath}
  \P( \{ | Y_{n-1}^3| > 1\} \cap B_{n-1}) \leq K \frac{(\log(n-1))^d}{(n-1)^{1+c/2}} \mathbf{1}_{d<1} + K \frac{(\log(n-1))^{d/ \eta}}{(n-1)^{1 / \eta +c/2\eta}} \mathbf{1}_{d \geq 1}.
\end{displaymath}

This finishes the proof. \Box

\section{ A Ferguson - Klass - LePage series representation of localisable processes in the $\sigma$-finite measure space case}\label{sigfinite}

When the space $E$ has infinite measure, one cannot use the representation above, since it is no longer possible to renormalize by $m(E)$. This is a major drawback, since typical applications we have in mind deal with processes defined on the real line, {\it i.e.} $E = \mathbb{R}$ and $m$ is the Lebesgue measure. However, in the $\sigma$-finite case, one may always perform a change of measure that allows to reduce to the finite case, as explained in \cite{ST}, proposition 3.11.3 (for specific examples of changes of measure, see section \ref{Examples}). In terms of localisability, this merely translates into adding a natural condition involving both the kernel and the change of measure:

\begin{theo}\label{msspfm2}
Let $(E,{\cal E},m)$ be a $\sigma$-finite measure space. Let $r : E \rightarrow \mathbb{R}_+$ be such that $\hat{m}(dx)=\frac{1}{r(x)}m(dx)$ is a probability measure. Let $\alpha$ be a $C^1$ function defined on $\bbbr$ and ranging in $[c,d] \subset (0,2)$. Let $b$ be a $C^1$ function defined on $\bbbr$. Let $f(t,u,.)$ be a family of functions such that, for all $(t,u) \in \bbbr^2$, $f(t,u,.) \in {\cal F}_{\alpha(u)}(E,{\cal E}, m)$. Let $(\Gamma_i)_{i \geq 1}$ be a sequence of arrival times of a Poisson process with unit arrival time, $(V_i)_{i \geq 1}$ be a sequence of i.i.d. random variables with distribution $\hat m $ on $E$, and $(\gamma_i)_{i \geq 1}$ be a sequence of i.i.d. random variables with distribution $P(\gamma_i=1)=P(\gamma_i=-1)=1/2$. Assume finally that the three sequences $(\Gamma_i)_{i \geq 1}$, $(V_i)_{i \geq 1}$, and $(\gamma_i)_{i \geq 1}$ are independent.
Consider the following random field:

\begin{equation}\label{msfm2}
X(t,u)= b(u)C^{1/\alpha(u)}_{\alpha(u)} \sum_{i=1}^{\infty} \gamma_i \Gamma_i^{-1/\alpha(u)} r(V_i)^{1/\alpha(u)}f(t,u,V_i),
\end{equation}
where $C_{\alpha} = \left( \int_{0}^{\infty} x^{-\alpha} \sin (x)dx \right)^{-1}$.
Assume that $X(t,u)$ (as a process in $t$) is localisable at $u$ with exponent $h \in (0,1)$ and local form $X'_u(t,u)$. Assume in addition that:

\begin{itemize}
    \item (Cs1) The family of functions $v \to f(t,v,x)$ is differentiable for all $(v,t)$ in a neighbourhood of $u$ and almost all $x$ in $E$. The derivatives of $f$ with respect to $v$ are denoted $f'_v$.
    \item (Cs2) There exists $\varep >0$ such that:
\begin{equation}\label{kercond1sf}
\sup_{t \in B(u,\varep)}  \int_E \sup_{w \in B(u,\varep)} (|f(t,w,x)|^{\alpha(w)}) \hspace{0.1cm} m(dx) < \infty.
\end{equation}
\item (Cs3) There exists $\varep >0$ such that:
\begin{equation}
\sup_{t \in B(u,\varep)}  \int_E \sup_{w \in B(u,\varep)} (|f'_u(t,w,x)|^{\alpha(w)}) \hspace{0.1cm} m(dx) < \infty.
\end{equation}

		\item (Cs4) There exists $\varep >0$ such that:
\begin{equation}\label{kercond2sf}
\sup_{t \in B(u,\varep)}  \int_E \sup_{w \in B(u,\varep)} \left[\left|f(t,w,x)\log|f(t,w,x)| \right|^{\alpha(w)}\right] \hspace{0.1cm} m(dx) < \infty.
\end{equation}
		\item (Cs5) There exists $\varep >0$ such that :
\begin{equation}\label{kercond3sf}
\sup_{t \in B(u,\varep)}  \int_E \sup_{w \in B(u,\varep)} \left[  \left|f(t,w,x) \log(r(x))  \right|^{\alpha(w)}  \right] \hspace{0.1cm} m(dx) < \infty.
\end{equation}
\end{itemize}

Then $Y(t) \equiv X(t,t)$ is localisable at $u$ with exponent $h$ and local form $Y'_u(t)=X'_u(t,u)$.

\end{theo}

\medskip

Remark: from (\ref{msfm2}), it may seem as though the process $Y$ depends on the particular change of measure used, {\it i.e.} the choice of a specific $r$. However, this is not case. More precisely, proposition \ref{fddsf} below shows that the finite dimensional distributions of $Y$ only depend on $m$.

\medskip

\noindent{\it Proof}

\medskip

We shall apply Theorem \ref{msspfm} to the function $g(t,w,x)=r(x)^{1 / \alpha(w)} f(t,w,x)$ on $(E, \mathcal{E}, \hat m )$.
\begin{itemize}
 \item By (Cs1), the family of functions $v \to f(t,v,x)$ is differentiable for all $(v,t)$ in a neighbourhood of $u$ and almost all $x$ in $E$ thus $v \to g(t,v,x)$  is differentiable too and (C1) holds.
\item Choose $\varep > 0 $ such that (Cs2) holds.
\begin{displaymath}
\sup_{w \in B(u,\varep)} (|g(t,w,x)|^{\alpha(w)}) = r(x) \sup_{w \in B(u,\varep)} (|f(t,w,x)|^{\alpha(w)}).
\end{displaymath}
One has
\begin{eqnarray*}
\int_\bbbr \sup_{w \in B(u,\varep)} (|g(t,w,x)|^{\alpha(w)}) \hspace{0.1cm} \hat m(dx) & =  & \int_\bbbr r(x) \sup_{w \in B(u,\varep)} (|f(t,w,x)|^{\alpha(w)}) \hspace{0.1cm} \hat m(dx)\\
 & = & \int_\bbbr \sup_{w \in B(u,\varep)} (|f(t,w,x)|^{\alpha(w)}) \hspace{0.1cm} m(dx)\\
\end{eqnarray*}
thus
\begin{eqnarray*}
\sup_{t \in B(u,\varep)} \int_\bbbr \sup_{w \in B(u,\varep)} (|g(t,w,x)|^{\alpha(w)}) \hspace{0.1cm} \hat m(dx) & =  & \sup_{t \in B(u,\varep)} \int_\bbbr \sup_{w \in B(u,\varep)} (|f(t,w,x)|^{\alpha(w)}) \hspace{0.1cm} m(dx)\\
\end{eqnarray*}
and (C2) holds.

\item Choose $\varep > 0 $ such that (Cs4) and (Cs5) hold.
\begin{displaymath}
\int_\bbbr \sup_{w \in B(u,\varep)} \left[\left|g(t,w,x)\log|g(t,w,x)|\right|^{\alpha(w)}\right] \hspace{0.1cm} \hat m(dx)
\end{displaymath}
\begin{displaymath}
 \leq \int_\bbbr r(x) \sup_{w \in B(u,\varep)} \left[\left|f(t,w,x)\log|r(x)^{1/ \alpha(w)}f(t,w,x)|\right|^{\alpha(w)}\right] \hspace{0.1cm} \hat m(dx)
\end{displaymath}
\begin{displaymath}
\leq \int_\bbbr \sup_{w \in B(u,\varep)} \left[\left|f(t,w,x)\log|r(x)^{1/ \alpha(w)}f(t,w,x)|\right|^{\alpha(w)}\right] \hspace{0.1cm} m(dx).
\end{displaymath}
Expanding the logarithm above and using the inequality $|a+b|^\delta \leq \max(1, 2^{\delta-1})(|a|^\delta + |b|^\delta)$, valid for all real numbers $a,b$ and all positive $\delta$, one sees that (C4) holds.

\item Choose $\varep > 0 $ such that (Cs3) and (Cs5) hold.
\begin{displaymath}
g'_u(t,w,x) = r(x)^{1/ \alpha(w)} \left(f'_u(t,w,x)-\frac{\alpha'(w)}{\alpha^2(w)}\log(r(x))f(t,w,x)\right)
\end{displaymath}
and
\begin{displaymath}
\int_\bbbr \sup_{w \in B(u,\varep)} (|g'_u(t,w,x)|^{\alpha(w)}) \hspace{0.1cm} \hat m(dx)
\end{displaymath}
\begin{displaymath}
\leq \int_\bbbr \sup_{w \in B(u,\varep)} \left[ \left|f'_u(t,w,x)-\frac{\alpha'(w)}{\alpha^2(w)}\log(r(x))f(t,w,x)\right|^{\alpha(w)} \right] \hspace{0.1cm} m(dx).
\end{displaymath}
The inequality $|a+b|^\delta \leq \max(1, 2^{\delta-1})(|a|^\delta + |b|^\delta)$ shows that (C3) holds.

Theorem \ref{msspfm} allows to conclude.
\end{itemize}

\Box

\section{Examples of localisable processes}\label{Examples}

In this section, we apply the results above and obtain some localisable
processes of interest. In particular, we consider ``multistable versions'' of several
classical processes. Similar multistable extensions were considered in
\cite{FLV}, to which the interested reader might refer for comparison.

We first recall some definitions. In the sequel, $M$ will denote
a symmetric $\alpha$-stable ($ 0 < \alpha < 2$)
random measure on $\bbbr$ with control measure Lebesgue measure ${\cal L}$.
We will write
$$L_{\alpha} (t) := \int_{0}^{t} M(dz)$$
for $\alpha$-stable L\'{e}vy motion.

The {\em log-fractional stable motion} is defined as
\begin{equation*}
\Lambda_\alpha(t)= \int_{-\infty}^{\infty} \left( \log(|t-x|) - \log(|x|)\right)
M(dx) \quad (t \in \bbbr).
\end{equation*}
This process is well-defined only for $\alpha \in (1,2]$
(the integrand does not belong to ${\cal F}_{\alpha}$ for $\alpha \leq 1$).
Both L\'evy motion and log-fractional stable motion
are $1/\alpha$-self-similar with stationary increments.

The following process is called \textit{linear fractional $\alpha$-stable motion}:
\begin{displaymath}
L_{\alpha,H,b^{+},b^{-}} (t) = \int_{-\infty}^{\infty} f_{\alpha,H}(b^{+},b^{-},t,x) M(dx)
\end{displaymath}
where $t \in \mathbb{R}$, $H \in (0,1)$, $b^{+},b^{-} \in \mathbb{R}$, and
\begin{align*}
f_{\alpha,H}(b^{+},b^{-},t,x) = b^{+} & \left( (t-x)_{+}^{H - 1/ \alpha} - (-x)_{+} ^{H - 1/ \alpha} \right) \nonumber
\\
& + b^{-} \left( (t-x)_{-} ^{H- 1/ \alpha} - (-x)_{-} ^{H- 1/ \alpha} \right). \label{lfsm}
\end{align*}
$L_{\alpha,H,b^{+},b^{-}}$ is again an sssi process. When $b^{+}=b^{-}=1$, this process
is called well-balanced linear fractional $\alpha$-stable motion and denoted $L_{\alpha,H}$.

Finally, for $\lambda>0$, the stationary process
\begin{equation*}
Y(t) = \int_{t}^{\infty} \exp(-\lambda (x-t))M(dx) \quad (t \in \bbbr)
\end{equation*}
is called reverse Ornstein-Uhlenbeck process.

The localisability of L\'{e}vy motion, log-fractional stable motion and
linear fractional $\alpha$-stable motion simply stems from the fact that
they are sssi. The localisability of the reverse Ornstein-Uhlenbeck process
is proved in \cite{FLGLV}.

We will now define multistable versions of these processes.

For the multistable L\'{e}vy motion, we give two versions: one is fitted to
the case where the time parameter varies in a compact interval $[0,T]$, and one
where it spans $\bbbr$.

\begin{theo} {\it (symmetric multistable L\'{e}vy motion, compact case)}.
Let $\alpha: [0,T] \to [c,d] \subset (1,2)$  and $b: [0,T] \to \bbbr^{+}$
be continuously differentiable. Let $(\Gamma_i)_{i \geq 1}$ be a sequence of arrival times of a Poisson process with unit arrival time,
$(V_i)_{i \geq 1}$ be a sequence of i.i.d. random variables with distribution $\hat{m}(dx)$, the uniform distribution on $[0,T]$,
and $(\gamma_i)_{i \geq 1}$ be a sequence of i.i.d. random variables with distribution $P(\gamma_i=1)=P(\gamma_i=-1)=1/2$.
Assume finally that the three sequences $(\Gamma_i)_{i \geq 1}$, $(V_i)_{i \geq 1}$, and $(\gamma_i)_{i \geq 1}$ are independent and
define
\begin{equation}
Y(t)  = b(t) C_{ \alpha(t)}^{1/ \alpha(t)} T^{1/ \alpha(t)}\sum_{i=1}^{+ \infty}
\gamma_i \Gamma_i^{-1/\alpha(t)} 1_{[0,t]}(V_i) \quad (t \in [0,T]).
\end{equation}

then  $Y$ is  $1/\alpha(u) $-localisable at any
$u \in (0,T)$,
with local form  $Y_{u}' = b(u)L_{\alpha(u)}$.

\end{theo}

\medskip

The proof is a simple application of Theorem \ref{msspfm}, and is omitted.

\medskip

\begin{theo} {\it (symmetric multistable L\'{e}vy motion, non-compact case)}.
Let $\alpha: \bbbr \to [c,d] \subset (1,2)$  and $b: \bbbr \to \bbbr^{+}$
be continuously differentiable. Let $(\Gamma_i)_{i \geq 1}$ be a sequence of arrival times of a Poisson process with unit arrival time,
$(V_i)_{i \geq 1}$ be a sequence of i.i.d. random variables with distribution $\hat{m}(dx) =\sum_{j=1}^{+\infty} 2^{-j} \mathbf{1}_{[j-1,j[}(x) dx$ on $\bbbr$,
and $(\gamma_i)_{i \geq 1}$ be a sequence of i.i.d. random variables with distribution $P(\gamma_i=1)=P(\gamma_i=-1)=1/2$.
Assume finally that the three sequences $(\Gamma_i)_{i \geq 1}$, $(V_i)_{i \geq 1}$, and $(\gamma_i)_{i \geq 1}$ are independent and
define
\begin{equation}
Y(t)  = b(t) C_{ \alpha(t)}^{1/ \alpha(t)} \sum_{i=1}^{+ \infty} \sum_{j =1}^{+\infty}
\gamma_i \Gamma_i^{-1/\alpha(t)} 2^{j/ \alpha(t)} 1_{[0,t] \cap [j-1,j[}(V_i) \quad (t \in \bbbr_+).
\end{equation}

then  $Y$ is  $1/\alpha(u) $-localisable at any
$u \in \bbbr_+$,
with local form  $Y_{u}' = b(u)L_{\alpha(u)}$.

\end{theo}

\medskip

\noindent{\it Proof}

\medskip

We apply Theorem \ref{msspfm2} with $m(dx)=dx$, $r(x)=\sum_{j=1}^{\infty} 2^j \mathbf{1}_{[j-1,j[}(x)$, $f(t,u,x)=\mathbf{1}_{[0,t]}(x)$ and the random field
\begin{displaymath}
 X(t,u)=b(u)C_{\alpha(u)}^{1/ \alpha(u)} \sum_{i,j=1}^{\infty} \gamma_i \Gamma_i^{-1/\alpha(u)} 2^{j/\alpha(u)} \mathbf{1}_{[0,t] \cap [j-1,j[}(V_i).
\end{displaymath}

$X(.,u)$ is the symmetrical $\alpha(u)$-L\'evy motion \cite{ST} and is thus $\frac{1}{\alpha(u)}$-localisable with local form $X'_u(.,u)=X(.,u)$.
\newline

\begin{itemize}
    \item (Cs1) The family of functions $v \to f(t,v,x)$ is differentiable for all $(v,t)$ in a neighbourhood of $u$ and almost all $x$ in $E$. The derivatives of $f$ with respect to $u$ vanish.
    \item (Cs2)
    \begin{displaymath}
  |f(t,w,x)|^{\alpha(w)} = \mathbf{1}_{[0,t]}(x)
    \end{displaymath}
thus
\begin{displaymath}
\int_\bbbr \sup_{w \in B(u,\varep)} (|f(t,w,x)|^{\alpha(w)}) \hspace{0.1cm} dx = t
\end{displaymath}
and (Cs2) holds.

\item (Cs3) $f'_u = 0$ so (Cs3) holds.

\item (Cs4) $ f(t,w,x)\log|f(t,w,x)| = 0$ so (Cs4) holds.

\item (Cs5)

\begin{eqnarray*}
 \left|f(t,w,x)\log(r(x))\right|^{\alpha(x)} &=& \sum_{j=1}^{+\infty} j^{\alpha(w)} \log(2)^{\alpha(w)} \mathbf{1}_{[0,t]\cap[j-1,j[}(x)\\
 & \leq & \log(2)^d \sum_{j=1}^{+\infty} j^d \mathbf{1}_{[0,t]\cap[j-1,j[}(x)\\
\end{eqnarray*}
thus
\begin{displaymath}
 \int_\bbbr \sup_{w \in B(u,\varep)} \left[\left|f(t,w,x)\log(r(x) )\right|^{\alpha(w)}\right] \hspace{0.1cm} dx \leq \log(2)^d \sum_{j=1}^{[t]+1} j^d
\end{displaymath}
and (Cs5) holds \Box

\end{itemize}
\Box

\begin{theo} {\it (Log-fractional multistable motion)}.
Let  $\alpha: \bbbr \to [c,d] \subset (1,2)$  and $b: \bbbr \to \bbbr^{+}$
be continuously differentiable. Let $(\Gamma_i)_{i \geq 1}$ be a sequence of arrival times of a Poisson process with unit arrival time,
$(V_i)_{i \geq 1}$ be a sequence of i.i.d. random variables with distribution $\hat{m}(dx) = \frac{3}{\pi^2}\sum_{j=1}^{+\infty} j^{-2} \mathbf{1}_{[-j,-j+1[ \cup [j-1,j[}(x) dx$ on $\bbbr$, and $(\gamma_i)_{i \geq 1}$ be a sequence of i.i.d. random variables with distribution $P(\gamma_i=1)=P(\gamma_i=-1)=1/2$. Assume finally that the three sequences $(\Gamma_i)_{i \geq 1}$, $(V_i)_{i \geq 1}$, and $(\gamma_i)_{i \geq 1}$ are independent and
define
\begin{equation}
Y(t)  = b(t) C_{ \alpha(t)}^{1/ \alpha(t)} \sum_{i=1}^{+ \infty} \sum_{j =1}^{+\infty}
\gamma_i \Gamma_i^{-1/\alpha(t)} (\log|t-V_i|-\log|V_i|)\frac{\pi^{2/\alpha(t)}}{3^{1/\alpha(t)}}j^{2/\alpha(t)} \mathbf{1}_{[-j,-j+1[ \cup [j-1,j[}(V_i)\quad (t \in \bbbr).
\end{equation}

then  $Y$ is  $1/\alpha(u) $-localisable at
any $u \in \bbbr$,
with $Y_{u}' = b(u)\Lambda_{\alpha(u)}$.

\end{theo}

\medskip

\noindent{\it Proof}

We apply Theorem \ref{msspfm2} with $m(dx)=dx$, $r(x)=\frac{\pi ^2}{3}\sum_{j=1}^{\infty} j^2 \mathbf{1}_{[-j,-j+1[ \cup [j-1,j[}(x)$, $f(t,u,x)=\log(|t-x|) - \log(|x|)$ and the random field
\begin{displaymath}
 X(t,u)=b(u)C_{\alpha(u)}^{1/ \alpha(u)} \sum_{i,j=1}^{\infty} \gamma_i \Gamma_i^{-1/\alpha(u)} (\log|t-V_i|-\log|V_i|)\frac{\pi^{2/ \alpha(u)}}{3^{1/ \alpha(u)}}j^{2/ \alpha(u)} \mathbf{1}_{[-j,-j+1[ \cup [j-1,j[}(V_i).
\end{displaymath}

$X(.,u)$ is the symmetrical $\alpha(u)$-Log-fractional motion. It is $\frac{1}{\alpha(u)}$-localisable with local form $X'_u(.,u)=b(u)\Lambda_{\alpha(u)}$ \cite{FLV}.
\newline

\begin{itemize}
    \item (Cs1) The family of functions $v \to f(t,v,x)$ is differentiable for all $(v,t)$ in a neighbourhood of $u$ and almost all $x$ in $E$. The derivatives of $f$ with respect to $u$ vanish.
    \item (Cs2) $\forall a > 1$, $\exists K_a > 0$ such that $\int_{\bbbr} |f(t,w,x)|^a dx \leq K_a|t|$ so
    \begin{eqnarray*}
  |f(t,w,x)|^{\alpha(w)} &=& |\log(|t-x|) - \log(|x|)|^{\alpha(w)}\\
   &=& |\log|1-\frac{t}{x}||^{\alpha(w)}\\
   & \leq & |\log|1-\frac{t}{x}||^{d} + |\log|1-\frac{t}{x}||^{c}\\
    \end{eqnarray*}
and
\begin{eqnarray*}
\int_\bbbr \sup_{w \in B(u,\varep)} (|f(t,w,x)|^{\alpha(w)}) \hspace{0.1cm} dx & \leq & (K_c+K_d)|t|\\
\end{eqnarray*}
thus (Cs2) holds.

\item (Cs3) $f'_u = 0$ so (Cs3) holds.

\item (Cs4)
\begin{eqnarray*}
|f(t,w,x)\log(|f(t,w,x)|)|^{\alpha(w)} & = & |f(t,w,x)\log(|f(t,w,x)|)|^{\alpha(w)}\mathbf{1}_{ \{|\log(|f(t,w,x)|)|\leq 1\}}\\
& &  +  |f(t,w,x)\log(|f(t,w,x)|)|^{\alpha(w)}\mathbf{1}_{\{|\log(|f(t,w,x)|)|> 1\}}\\
& \leq & |f(t,w,x)|^{\alpha(w)} + |f(t,w,x)\log(|f(t,w,x)|)|^{\alpha(w)} \mathbf{1}_{\{|f(t,w,x)|>e\}}\\
& &  + |f(t,w,x)\log(|f(t,w,x)|)|^{\alpha(w)}\mathbf{1}_{\{|f(t,w,x)| < \frac{1}{e}\}}.\\
\end{eqnarray*}
We shall bound each of the three terms that are added up in the right hand side
of the above inequality.
For the first term,
\begin{displaymath}
|f(t,w,x)|^{\alpha(w)} \leq \sup_{w \in B(u,\varep)}|f(t,w,x)|^{\alpha(w)}.
\end{displaymath}

For the second term, fix $K>0$, $\epsilon > 0$ such that $\forall x > e, |x\log(|x|)| \leq K|x|^{1+\epsilon}$.
\begin{eqnarray*}
|f(t,w,x)\log(|f(t,w,x)|)|^{\alpha(w)} \mathbf{1}_{\{|f(t,w,x)|>e\}} & \leq & K|f(t,w,x)|^{d(1+\epsilon)}.\\
\end{eqnarray*}

For the third term,
fix $K_1 < K_2 < 0$ and $K_4 > K_3 >0$ such that

\begin{displaymath}
\mathbf{1}_{\{|f(t,w,x)| < \frac{1}{e}\}} \leq \mathbf{1}_{]-\infty,K_1|t|[}(x) + \mathbf{1}_{]K_2|t|, K_3|t|[}(x) + \mathbf{1}_{]K_4|t|, +\infty[}(x),
\end{displaymath}
then

\begin{eqnarray*}
|f(t,w,x)\log(|f(t,w,x)|)|^{\alpha(w)}\mathbf{1}_{\{|f(t,w,x)| < \frac{1}{e}\}} & \leq & |f(t,w,x)|^c |\log(|f(t,w,x)|)|^d \mathbf{1}_{]-\infty,K_1|t|[}(x)\\
& & + |f(t,w,x)|^c |\log(|f(t,w,x)|)|^d \mathbf{1}_{]K_2|t|, K_3|t|[}(x)\\
& & + |f(t,w,x)|^c |\log(|f(t,w,x)|)|^d\mathbf{1}_{]K_4|t|, +\infty[}(x).\\
\end{eqnarray*}

The function $x \mapsto |x|^c|\log(|x|)|^d $ is bounded for $x < \frac{1}{e} $. With $M$ denoting an upper bound of this function, one has

\begin{eqnarray*}
\int_{K_2|t|}^{K_3|t|}|f(t,w,x)|^c |\log(|f(t,w,x)|)|^d  dx & \leq & M (K_3-K_2)|t|.\\
\end{eqnarray*}

With $u = 1 - \frac{t}{x}$, we obtain
\begin{eqnarray*}
\int_{K_4|t|}^{+\infty} |f(t,w,x)|^c |\log(|f(t,w,x)|)|^d dx & \leq  &  K'|t| + |t|\int_{1 - \frac{1}{K''}}^{+\infty} \frac{|\log|u||^c|\log|\log|u|||^d}{(1-u)^2} du\\
& \leq & O(|t|),\\
\end{eqnarray*}
where $K'$ and $K''$ are numbers verifying $K''>1$ and $K'>0$.

For the same reason, $ \int_{-\infty}^{K_1|t|} |f(t,w,x)|^c |\log(|f(t,w,x)|)|^d dx \leq K|t|$.
We conclude that
\begin{displaymath}
\sup_{t \in B(u,\varep)}  \int_\bbbr \sup_{w \in B(u,\varep)} \left[\left|f(t,w,x)\log(|f(t,w,x)|)\right|^{\alpha(w)}\right] \hspace{0.1cm} dx < \infty.
\end{displaymath}

\item (Cs5)
\begin{eqnarray*}
 |f(t,w,x)\log(r(x))|^{\alpha(w)}  & \leq & K_1|f(t,w,x)\log(\frac{\pi^2}{3})|^{\alpha(w)}\\
& & + K_2 \sum_{j=1}^{+\infty} |f(t,w,x)|^{\alpha(w)} (\log(j))^d \mathbf{1}_{[-j,-j+1[ \cup [j-1,j[}(x).\\
\end{eqnarray*}
For $j$ large enough ($j > j^*$), $|f(t,w,x)|^{\alpha(w)}\mathbf{1}_{[-j,-j+1[ \cup [j-1,j[}(x) \leq K_5 \frac{|t|^c}{|x|^c} \mathbf{1}_{[-j,-j+1[ \cup [j-1,j[}(x)$). Thus

\begin{eqnarray*}
 |f(t,w,x)\log(r(x))|^{\alpha(w)}  & \leq & K_6|f(t,w,x)|^{\alpha(w)} + K_7 \sum_{j=j^*}^{+\infty} (\log(j))^d \frac{|t|^c}{|x|^c}\mathbf{1}_{[-j,-j+1[ \cup [j-1,j[}(x).\\
\end{eqnarray*}

To conclude, note that
\begin{eqnarray*}
\int_{\bbbr} (\log(j))^d \frac{1}{|x|^c}\mathbf{1}_{[-j,-j+1[ \cup [j-1,j[}(x) dx & = & 2 (\log(j))^d \int_{j-1}^{j} \frac{dx}{|x|^c}\\
& \sim & 2 (c-1) \frac{(\log(j))^d}{j^c}\\
\end{eqnarray*}

\Box

\end{itemize}

\begin{theo} {\it (Linear multistable multifractional motion)}.
Let  $b : \bbbr \to \bbbr^{+}$
, $\alpha:\bbbr \to [c,d] \subset (0,2)$ and
$h:\bbbr \to (0,1)$ be continuously differentiable.
 Let $(\Gamma_i)_{i \geq 1}$ be a sequence of arrival times of a Poisson process with unit arrival time,
$(V_i)_{i \geq 1}$ be a sequence of i.i.d. random variables with distribution $\hat{m}(dx) = \frac{3}{\pi^2}\sum_{j=1}^{+\infty} j^{-2} \mathbf{1}_{[-j,-j+1[ \cup [j-1,j[}(x) dx$ on $\bbbr$, and $(\gamma_i)_{i \geq 1}$ be a sequence of i.i.d. random variables with distribution $P(\gamma_i=1)=P(\gamma_i=-1)=1/2$. Assume finally that the three sequences $(\Gamma_i)_{i \geq 1}$, $(V_i)_{i \geq 1}$, and $(\gamma_i)_{i \geq 1}$ are independent and
define for $t \in \bbbr$
\begin{equation}
Y(t)  = b(t) C_{ \alpha(t)}^{1/ \alpha(t)} \sum_{i,j=1}^{+ \infty}
\gamma_i \Gamma_i^{-1/\alpha(t)} (|t-V_i|^{h(t)-1/ \alpha(t)}-|V_i|^{h(t)-1/\alpha(t)})(\frac{\pi^2 j^2}{3})^{1/ \alpha(t)} \mathbf{1}_{[-j,-j+1[ \cup [j-1,j[}(V_i).
\end{equation}

The process  $Y$ is $h(u)$-localisable
at all $u \in \bbbr$, with $Y_{u}' = b(u)L_{\alpha(u),h(u)}$ (the well balanced linear fractional stable motion).
\end{theo}

\noindent{\it Proof}

We apply Theorem \ref{msspfm2} with $m(dx)=dx$, $r(x)=\frac{\pi ^2}{3}\sum_{j=1}^{\infty} j^2 \mathbf{1}_{[-j,-j+1[ \cup [j-1,j[}(x)$, $f(t,u,x)=|t-x|^{h(u)-1/ \alpha(u)}-|x|^{h(u)-1/\alpha(u)}$ and the random field
\begin{displaymath}
 X(t,u)=b(u)C_{\alpha(u)}^{1/ \alpha(u)} (\frac{\pi^2 j^2}{3})^{1/ \alpha(u)}\sum_{i,j=1}^{\infty} \gamma_i \Gamma_i^{-1/\alpha(u)} (|t-V_i|^{h(u)-1/ \alpha(u)}-|V_i|^{h(u)-1/\alpha(u)}) \mathbf{1}_{[-j,-j+1[ \cup [j-1,j[}(V_i).
\end{displaymath}

$X(.,u)$ is the ($\alpha(u),h(u)$)-well balanced linear fractional stable motion and it is $\frac{1}{\alpha(u)}$-localisable with local form $X'_u(.,u)=b(u)L_{\alpha(u),h(u)}$ \cite{FLV}.
\newline

\begin{itemize}
    \item (Cs1) The family of functions $v \to f(t,v,x)$ is differentiable for all $(v,t)$ in a neighbourhood of $u$ and almost all $x$ in $E$. The derivatives of $f$ with respect to $u$ read:
\begin{displaymath}
f'_u(t,w,x)=(h'(w)+\frac{\alpha'(w)}{\alpha^2(w)})\left[ (\log|t-x|)|t-x|^{h(w)-1/\alpha(w)} - (\log|x|)|x|^{h(w)-1/\alpha(w)}\right].
\end{displaymath}
    \item (Cs2) In \cite{FLV}, it is shown that, given $u \in \bbbr$, one may choose $\varep>0$ small enough and numbers $a,b,h_{-},h_{+}$ with
$0<a <\alpha(w) < b<2 $, $0<h_{-} < h(w) < h_{+} <1$ and
$\frac{1}{a}-\frac{1}{b}< h_{-}<h_{+}<1-(\frac{1}{a}-\frac{1}{b})$
such that, for all $t$ and $w$ in $U:=(u-\varep, u+\varep)$ and all real $x$:
\begin{equation}\label{inegk}
|f(t,w,x)|,
|f'_{u}(t,w,x)|
\leq  k_{1}(t,x)
\end{equation}

 where
\begin{equation}
k_{1}(t,x) = \left\{
\begin{array}{ll}
    c_{1}\max\{1, |t-x|^{h_{-}-1/a} +|x|^{h_{-}-1/a}\}
    & ( |x| \leq 1 + 2 \max_{t \in U}|t|) \\
    c_{2}|x|^{h_{+}-1/b-1}
     & ( |x| > 1 + 2 \max_{t \in U}|t|)
\end{array}
\right.
\end{equation}
for appropriately chosen constants $c_{1}$ and $c_{2}$. The conditions on $a,b,h_{-},h_{+}$ entail that
$\sup_{t \in U}\| k_{1}(t,\cdot)\|_{a,b}< \infty$ and (Cs2) hold.

    \item (Cs3) is obtained with (\ref{inegk}) for the same reason as in (Cs2).

    \item (Cs4)

    \begin{eqnarray*}
|f(t,w,x)\log(|f(t,w,x)|)|^{\alpha(w)} & \leq & |f(t,w,x)|^{\alpha(w)} + |f(t,w,x)\log(|f(t,w,x)|)|^{\alpha(w)} \mathbf{1}_{\{|f(t,w,x)|>e\}}\\
& &  + |f(t,w,x)\log(|f(t,w,x)|)|^{\alpha(w)}\mathbf{1}_{\{|f(t,w,x)| < \frac{1}{e}\}}.\\
\end{eqnarray*}
Since
    \begin{equation}
|f(t,w,x)| \leq  k_{1}(t,x)
\end{equation}
one gets

\begin{eqnarray*}
  |f(t,w,x)\log(|f(t,w,x)|)|\mathbf{1}_{\{|f(t,w,x)|>e\}} & \leq & k_{1}(t,x) \log(k_{1}(t,x))\mathbf{1}_{\{|f(t,w,x)|>e\}}\\
  & \leq  & |k_{1}(t,x) \log(k_{1}(t,x))|\\
 |f(t,w,x)\log(|f(t,w,x)|)|^{\alpha(w)} \mathbf{1}_{\{|f(t,w,x)|>e\}} & \leq &  |k_{1}(t,x)\log(k_{1}(t,x))|^{\alpha(w)}.\\
\end{eqnarray*}

Fix $\eta >0$ such that $1 < \eta < a + \frac{a}{b} - a h_+$ and $\lambda > 0$ such that $\frac{1}{\eta} < \lambda < 1$.
\begin{eqnarray*}
|f(t,w,x)\log(|f(t,w,x)|)|^{\alpha(w)}\mathbf{1}_{\{|f(t,w,x)| < \frac{1}{e}\}} & \leq & K |f(t,w,x)|^{\lambda \alpha(w)}\\
& \leq & K |k_{1}(t,x)|^{\lambda \alpha(w)}\\
\end{eqnarray*}
and thus (Cs4) holds.

    \item (Cs5) For $j$ large enough ($j > j^*$),
    \begin{eqnarray*}
 |f(t,w,x)\log(r(x))|^{\alpha(w)}  & \leq & K_1|k_1(t,x)|^{\alpha(w)}\\
& & + K_2 \sum_{j=j^*}^{+\infty} |f(t,w,x)|^{\alpha(w)} (\log(j))^d \mathbf{1}_{[-j,-j+1[ \cup [j-1,j[}(x).\\
\end{eqnarray*}
\begin{eqnarray*}
 |f(t,w,x)|^{\alpha(w)}\mathbf{1}_{[-j,-j+1[ \cup [j-1,j[}(x) &\leq& K_3 \frac{1}{|x|^{a(1
 +1/b-h_+)}} \mathbf{1}_{[-j,-j+1[ \cup [j-1,j[}(x).
\end{eqnarray*}
Thus

\begin{eqnarray*}
 |f(t,w,x)\log(r(x))|^{\alpha(w)}  & \leq & K_1|k_1(t,x)|^{\alpha(w)} + K_4 \sum_{j=j^*}^{+\infty} \frac{(\log(j))^d}{|x|^{a(1
 +1/b-h_+)}}\mathbf{1}_{[-j,-j+1[ \cup [j-1,j[}(x).\\
\end{eqnarray*}

To conclude, note that
\begin{eqnarray*}
\int_{\bbbr} \frac{(\log(j))^d}{|x|^{a(1
 +1/b-h_+)}}\mathbf{1}_{[-j,-j+1[ \cup [j-1,j[}(x) dx & = & 2 (\log(j))^d \int_{j-1}^{j} \frac{dx}{|x|^{a(1+1/b-h_+)}}\\
& \sim & 2 ({a(1+1/b-h_+)}-1) \frac{(\log(j))^d}{j^{a(1+1/b-h_+)}}\\
\end{eqnarray*}

\Box

\end{itemize}

\begin{theo} {\it (Multistable reverse Ornstein-Uhlenbeck process)}.
Let $\lambda > 0$, $\alpha: \bbbr \to [c,d] \subset (1,2)$  and $b: \bbbr \to \bbbr^{+}$
be continuously differentiable. Let $(\Gamma_i)_{i \geq 1}$ be a sequence of arrival times of a Poisson process with unit arrival time, $(V_i)_{i \geq 1}$ be a sequence of i.i.d. random variables with distribution $\hat{m}(dx) = \sum_{j=1}^{+\infty} 2^{-j-1} \mathbf{1}_{[-j,-j+1[ \cup [j-1,j[}(x) dx$ on $\bbbr$, and $(\gamma_i)_{i \geq 1}$ be a sequence of i.i.d. random variables with distribution $P(\gamma_i=1)=P(\gamma_i=-1)=1/2$. Assume finally that the three sequences $(\Gamma_i)_{i \geq 1}$, $(V_i)_{i \geq 1}$, and $(\gamma_i)_{i \geq 1}$ are independent and
define
\begin{equation}
Y(t)  = b(t) C_{ \alpha(t)}^{1/ \alpha(t)} \sum_{i=1}^{+ \infty} \sum_{j =1}^{+\infty}
\gamma_i \Gamma_i^{-1/\alpha(t)} 2^{(j+1)/ \alpha(t)} e^{-\lambda(V_i-t)} \mathbf{1}_{[t,+\infty) \cap ([-j,-j+1)\cup [j-1,j))}(V_i) \quad (t \in \bbbr).
\end{equation}

Then  $Y$ is  $1/\alpha(u) $-localisable at any
$u \in \bbbr $,
with local form $Y_{u}' = b(u)L_{\alpha(u)}$.

\end{theo}

\medskip

\noindent{\it Proof}

We apply Theorem \ref{msspfm2} with $m(dx)=dx$, $r(x)=\sum_{j=1}^{\infty} 2^{j+1} \mathbf{1}_{[-j,-j+1) \cup [j-1,j)}(x)$, $f(t,u,x)=e^{-\lambda(x-t)}\mathbf{1}_{[t,+\infty)}(x)$ and the random field
\begin{displaymath}
 X(t,u)=b(u)C_{\alpha(u)}^{1/ \alpha(u)} \sum_{i,j=1}^{\infty} \gamma_i \Gamma_i^{-1/\alpha(u)} 2^{(j+1)/\alpha(u)}e^{-\lambda(V_i-t)} \mathbf{1}_{[t,+\infty) \cap ([-j,-j+1) \cup [j-1,j))}(V_i).
\end{displaymath}

$X(.,u)$ is the symmetrical $\alpha(u)$-reverse Ornstein-Uhlenbeck process and is $\frac{1}{\alpha(u)}$-localisable with local form $X'_u(.,u)=b(u)L_{\alpha(u)}$ \cite{FLGLV}.
\newline

\begin{itemize}
    \item (Cs1) The family of functions $v \to f(t,v,x)$ is differentiable for all $(v,t)$ in a neighbourhood of $u$ and almost all $x$ in $E$. The derivatives of $f$ with respect to $u$ vanish.
    \item (Cs2)
    \begin{displaymath}
  |f(t,w,x)|^{\alpha(w)} = e^{-\lambda \alpha(w)(x-t)}\mathbf{1}_{[t,+\infty[}(x)
    \end{displaymath}
and
\begin{eqnarray*}
\int_\bbbr \sup_{w \in B(u,\varep)} (|f(t,w,x)|^{\alpha(w)}) \hspace{0.1cm} dx & \leq & \int_{t}^{+\infty} e^{-\lambda c (x-t)} dx\\
& \leq & \frac{1}{\lambda c}\\
\end{eqnarray*}
thus (Cs2) holds.
\item (Cs3) $f'_u = 0$ so (Cs3) holds.

\item (Cs4)

\begin{displaymath}
|f(t,w,x)\log(|f(t,w,x)|)|^{\alpha(w)} = \lambda^{\alpha(w)} (x-t)^{\alpha(w)} e^{-\lambda \alpha(w) (x-t)} \mathbf{1}_{[t,+\infty[}(x)
\end{displaymath}
as a consequence
\begin{eqnarray*}
 \int_\bbbr \sup_{w \in B(u,\varep)} \left[\left|f(t,w,x)\log(|f(t,w,x)|)\right|^{\alpha(w)}\right] \hspace{0.1cm} dx & \leq & \int_{0}^{+\infty} \lambda^d u^d e^{-\lambda c u} \hspace{0.1cm} du\\
 & < & +\infty.\\
\end{eqnarray*}

\item (Cs5)
\begin{displaymath}
 |f(t,w,x)\log(r(x))|^{\alpha(w)} = \sum_{j=1}^{+\infty} (j+1)^{\alpha(w)} \log(2)^{\alpha(w)} e^{-\lambda \alpha(w) (x-t)} \mathbf{1}_{[t,+\infty[ \cap ([-j,-j+1[ \cup [j-1,j[)} (x).
\end{displaymath}
Fix $j^*$ large enough such that for all $j > j^*$, $\mathbf{1}_{[t,+\infty[ \cap ([-j,-j+1[ \cup [j-1,j[)} (x) = \mathbf{1}_{[j-1,j[} (x)$. Then
\begin{eqnarray*}
 \int_\bbbr \sup_{w \in B(u,\varep)} \left[\left|f(t,w,x)\log(r(x))\right|^{\alpha(w)}\right] \hspace{0.1cm} dx & \leq &\sum_{j=1}^{j^*} \frac{(j+1)^d \log(2)^d}{\lambda c}\\
 \end{eqnarray*}
\begin{displaymath}
 + \sum_{j=j^* +1}^{+\infty} (j+1)^d \log(2)^d \int_{j-1}^{j} e^{-\lambda c (x-t)} \hspace{0.1cm} dx
\end{displaymath}
\begin{eqnarray*}
& \leq & \sum_{j=1}^{j^*} \frac{(j+1)^d \log(2)^d}{\lambda c} + \log(2)^d e^{\lambda c t} ( e^{\lambda c}-1) \sum_{j=j^* +1}^{+\infty} (j+1)^d e^{-\lambda c j}\\
 \end{eqnarray*}

\Box

\end{itemize}

\section{Finite dimensional distributions}\label{fdd}

In this section, we compute the finite dimensional distributions of the family of processes
defined in theorem \ref{msspfm2}, and compare the results with the ones in \cite{FLV}.

\begin{prop}\label{fddsf}
With notations as above, let $\{ X(t,u), t,u \in \bbbr \}$  be as in $(\ref{msfm2})$ and $Y(t) \equiv X(t,t)$. The finite dimensional distributions of the process $Y$ are equal to
\begin{displaymath}
\E \left( e^{i \sum_{j=1}^{m}\limits \theta_j Y(t_j)} \right) = \exp \left( -2 \int_{E} \int_{0}^{+ \infty} \sin^2( \sum_{j=1}^{m} \theta_j b(t_j) \frac{C_{\alpha(t_j)}^{1/\alpha(t_j)}}{2y^{1/\alpha(t_j)}} f(t_j,t_j,x) )\hspace{0.1cm} dy \hspace{0.1cm}m(dx) \right)
\end{displaymath}
for $ m \in \mathbb{N}, \boldsymbol{\theta} = (\theta_1, \ldots, \theta_m) \in \bbbr^m, {\bf t } = (t_1, \ldots , t_m) \in \bbbr^m$.
\end{prop}

\medskip

\noindent{\it Proof.}
Let $m \in \mathbb{N}$ and write $\phi_t(\theta) = \E \left(e^{i \sum_{j=1}^{m}\limits \theta_j Y(t_j)}\right)$. We proceed as in \cite{ST}, proposition 1.4.2. Let $\{U_i\}_{i \in \mathbb{N}}$ be an i.i.d sequence of uniform random variables on $(0,1)$, independent of the sequences $\{\gamma_i\}$ and $\{V_i\}$, and $g(t,u,x)=b(u) C_{\alpha(u)}^{1/ \alpha(u)} r(x)^{1/ \alpha(u)} f(t,u,x)$. For all $n \in \mathbb{N}$,

\begin{equation}\label{eqdist}
\sum_{j=1}^{m} \theta_j n^{-1/ \alpha(t_j)} \sum_{k=1}^{n} \gamma_k U_k^{-1/ \alpha(t_j)} g(t_j,t_j,V_k) \ed \sum_{j=1}^{m} \theta_j \left( \frac{\Gamma_{n+1}}{n}\right)^{1/ \alpha(t_j)} \sum_{k=1}^{n} \gamma_k \Gamma_k^{-1/ \alpha(t_j)} g(t_j,t_j,V_k).
\end{equation}
The right-hand side of (\ref{eqdist}) converges almost surely
to $\sum_{j=1}^{m}\limits \theta_j Y(t_j)$ when $n$ tends to infinity and thus
\begin{displaymath}
\phi_t(\theta) = \lim_{n \rightarrow +\infty} \E \left(e^{i \sum_{j=1}^{m}\limits \theta_j n^{-1/ \alpha(t_j)} \sum_{k=1}^{n}\limits  \gamma_k U_k^{-1/ \alpha(t_j)} g(t_j,t_j,V_k)}\right).
\end{displaymath}
Set $\phi_t^n(\theta) = \E \left(e^{i \sum_{j=1}^{m}\limits \theta_j n^{-1/ \alpha(t_j)} \sum_{k=1}^{n}\limits  \gamma_k U_k^{-1/ \alpha(t_j)} g(t_j,t_j,V_k)}\right)$. This
function may be written as:
\begin{displaymath}
\phi_t^n(\theta) = \E \left(\prod_{k=1}^{n} e^{i \gamma_k \sum_{j=1}^{m}\limits \theta_j n^{-1/ \alpha(t_j)}  U_k^{-1/ \alpha(t_j)} g(t_j,t_j,V_k)}\right).
\end{displaymath}
All the sequences $\{ \gamma_k\}$, $\{U_k \}$, $\{V_k \}$ are i.i.d. As a consequence,
\begin{displaymath}
\phi_t^n(\theta) = \left( \E \left(e^{i \gamma_1 \sum_{j=1}^{m}\limits \theta_j n^{-1/ \alpha(t_j)}  U_1^{-1/ \alpha(t_j)} g(t_j,t_j,V_1)}\right)\right)^n.
\end{displaymath}

We compute now the expectation using conditioning and independence of the sequences $\{ \gamma_k\}$, $\{U_k \}$ and $\{V_k \}$.

\begin{eqnarray*}
\E \left(e^{i \gamma_1 \sum_{j=1}^{m}\limits \theta_j n^{-1/ \alpha(t_j)}  U_1^{-1/ \alpha(t_j)} g(t_j,t_j,V_1)}\right) & = & \E\left( \E\left( e^{i \gamma_1 \sum_{j=1}^{m}\limits \theta_j n^{-1/ \alpha(t_j)}  U_1^{-1/ \alpha(t_j)} g(t_j,t_j,V_1)} | U_1,V_1\right)\right)\\
& = & \E\left( \cos(\sum_{j=1}^{m}\limits \theta_j n^{-1/ \alpha(t_j)}  U_1^{-1/ \alpha(t_j)} g(t_j,t_j,V_1))\right)\\
& = & \E \left(\frac{1}{n} \int_0^n \cos (\sum_{j=1}^{m}\limits \theta_j y^{-1/ \alpha(t_j)} g(t_j,t_j,V_1)) \hspace{0.1cm} dy\right)\\
& = & 1-\frac{2}{n} \int_0^n \E \left(\sin^2(\sum_{j=1}^{m}\limits \frac{\theta_j}{2} y^{-1/ \alpha(t_j)} g(t_j,t_j,V_1))\right) \hspace*{0.1cm} dy.\\
\end{eqnarray*}
The function $\sin^2$ is positive and thus, when $n$ tends to $+\infty$,
\begin{displaymath}
 \int_0^n \E\left(\sin^2(\sum_{j=1}^{m}\limits \frac{\theta_j}{2} y^{-1/ \alpha(t_j)} g(t_j,t_j,V_1))\right) \hspace*{0.1cm} dy \rightarrow  \int_0^{+\infty} \E \left(\sin^2(\sum_{j=1}^{m}\limits \frac{\theta_j}{2} y^{-1/ \alpha(t_j)} g(t_j,t_j,V_1))\right) \hspace*{0.1cm} dy.
\end{displaymath}
To conclude, note that
\begin{displaymath}
\E \left(\sin^2(\sum_{j=1}^{m}\limits \frac{\theta_j}{2} y^{-1/ \alpha(t_j)} g(t_j,t_j,V_1))\right)= \int_E \sin^2(\sum_{j=1}^{m}\limits \frac{\theta_j}{2} y^{-1/ \alpha(t_j)} g(t_j,t_j,x)) \hspace{0.1cm} \hat{m}(dx).
\end{displaymath}
\Box

Comparing with proposition 8.2, Theorems 9.3, 9.4, 9.5 and 9.6 in \cite{FLV}, it is easy to prove the following corollary, which shows that the approach based on the series representation and the
one based on sums over Poisson processes yield essentially the same processes:

\begin{cor}
The linear multistable multifractional motion, multistable L\'evy motion, log-fractional multistable
motion and multistable reverse Ornstein-Uhlenbeck process defined in section \ref{Examples} have the
same finite dimensional distributions as the corresponding processes considered in \cite{FLV}.
\end{cor}

\section{Numerical experiments}\label{numexp}

We display in this section graphs of synthesized paths of some of the processes
defined above. The idea is just to picture how multistability translates on the
behaviour of random trajectories, and, in the case of linear multistable multifractional
motion, to visualize the effect of both a varying $H$ and a varying $\alpha$,
these two parameters corresponding to two different notions of irregularity.
The synthesis method is described in \cite{FLGLV}. Theoretical results concerning
the convergence of this method will be presented elsewhere.

\medskip

\begin{figure}[ht]
\vspace{-0.15cm}
\includegraphics[height=18.9cm]{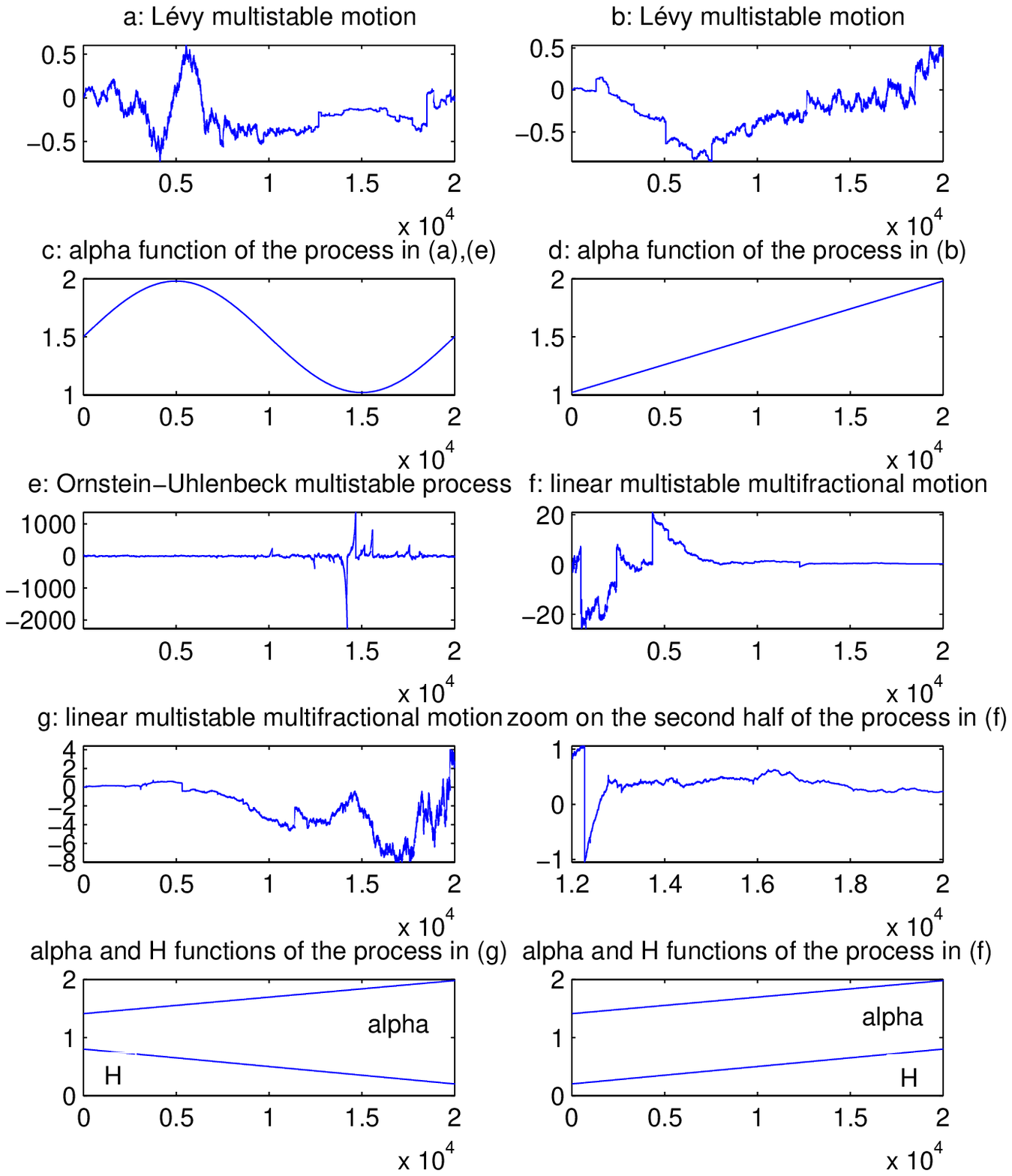}
\caption{Paths of multistable processes. First line: Levy multistable
motions with sine (a) and linear (b) $\alpha$ function. Second line:
(c) $\alpha$ function for the process in (a), (d) $\alpha$ function
for the process in (b). Third line: (e) multistable Ornstein-Uhlenbeck process
with $\alpha$ function displayed in (c), and (f) linear multistable
multifractional motion with linear increasing $\alpha$ and $H$ functions.
Fourth line: (g) linear multistable multifractional motion with linear
increasing $\alpha$ function and linear decreasing $H$ function, and zoom
on the second part of the process in (f). Last line: $\alpha$ and $H$
functions for the process in (g) (left), $\alpha$ and $H$
functions for the process in (f) (right).}
\label{tout}
\end{figure}

The two graphs on the first line of Figure \ref{tout} ((a) and (b)) display multistable L\'evy motions, with respectively
$\alpha$ increasing linearly from 1.02 to 1.98 (shown in (c)) and $\alpha$ a sine function ranging in
the same interval (shown in (d)). The graph (e) displays an Ornstein-Uhlenbeck multistable process
with same sine $\alpha$ function. A linear multistable multifractional
motion with linearly increasing $\alpha$ and $H$ functions is shown in (f). $H$ increases from
0.2 to 0.8 and $\alpha$ from 1.41 to 1.98 (these two functions are displayed on the right part of
the bottom line). The graph in (g) is again a linear multistable multifractional
motion, but with linearly increasing $\alpha$ and linearly decreasing $H$. $H$ decreases from
0.8 to 0.2 and $\alpha$ increases from 1.41 to 1.98 (these two functions are displayed on the left
part of the bottom line). Finally, a zoom on the second half of the process in (f) is shown, that
allows to appreciate how the graph becomes smoother as $H$ increases.

\medskip

In all these graphs, one clearly sees how the variations of $\alpha$ translates in terms of
the ``intensity'' of jumps. Additionally, in the case of linear multistable multifractional
motions, the interplay between the smoothness governed by $H$ and the jumps tuned by $\alpha$
indicate that such processes may prove useful in various applications such as finance or biomedical
signal modeling.

\end{document}